%% file: main_document.tex
\documentclass[11pt]{article}

    \newcommand{\blind}{0}
    
    \makeatletter
    \renewcommand\section{\@startsection {section}{1}{\z@}
                                       {-3.5ex \@plus -1ex \@minus -.2ex}
                                       {2.3ex \@plus.2ex}
                                       {\normalfont\fontfamily{phv}\fontsize{16}{19}\bfseries}}
    \renewcommand\subsection{\@startsection{subsection}{2}{\z@}
                                         {-3.25ex\@plus -1ex \@minus -.2ex}
                                         {1.5ex \@plus .2ex}
                                         {\normalfont\fontfamily{phv}\fontsize{14}{17}\bfseries}}
    \renewcommand\subsubsection{\@startsection{subsubsection}{3}{\z@}
                                        {-3.25ex\@plus -1ex \@minus -.2ex}
                                         {1.5ex \@plus .2ex}
                                         {\normalfont\normalsize\fontfamily{phv}\fontsize{14}{17}\selectfont}}
    \makeatother
   
	\usepackage{amsmath}
	\usepackage{graphicx}
	\usepackage{enumerate}
	\usepackage{xcolor}
	\usepackage{natbib} 

        \usepackage{url}
	
        \usepackage{amsfonts, amsthm, latexsym, amssymb}
        \usepackage{lineno}
        
        \newtheorem{theorem}{Theorem}[]
        
        \newtheorem{corollary}{Corollary}[]
        \usepackage{enumitem}
        \usepackage{booktabs,caption}
        \usepackage{wrapfig}
        \usepackage{float}
        \usepackage{fancyhdr}                                
        \usepackage[flushleft]{threeparttable}
        \usepackage[ruled]{algorithm2e}
        \SetAlgoNlRelativeSize{-1}
        \SetAlgoSkip{smallskip}
        \usepackage{mathtools}
        \usepackage{xcolor}
        \usepackage{arydshln}
        \usepackage{caption}
        \usepackage[resetlabels]{multibib}
        \newcites{Supp}{References}

        \definecolor{mycolor}{RGB}{0, 0, 0}

\begin{document}

		\def\spacingset#1{\renewcommand{\baselinestretch}
			{#1}\small\normalsize} \spacingset{1}
	
		\if0\blind
		{
			\title{\bf \emph{Effects of Geopolitical Strain on Global Pharmaceutical Supply Chain Design and Drug Shortages}} 
			\author{Martha L. Sabogal De La Pava $^a$ and Emily L. Tucker $^{a,b}$ \\
			$^a$ Industrial Engineering Department, Clemson University, Clemson, United States \\
                $^b$ School of Health Research, Clemson University, Clemson, United States }
			\date{}
			\maketitle
		} \fi
		
		\if1\blind
		{

            \title{\bf \emph{Effects of Geopolitical Strain on Global Pharmaceutical Supply Chain Design and Drug Shortages}} 
			\author{Author information is purposely removed for double-blind review}
			
        \bigskip
			\bigskip
			\bigskip
			\begin{center}
				{\LARGE\bf \emph{Effects of Geopolitical Strain on Global Pharmaceutical Supply Chain Design and Drug Shortages}} 
			\end{center}
			\medskip
		} \fi
		\bigskip

	\begin{abstract}
        Emerging geopolitical risks have begun to threaten global supply chains, including those that produce life-saving drugs. Export bans may prevent a company from shipping products internationally, and it is unclear how these new dynamics may affect company plans and persistent, worldwide drug shortages. To address these questions, we present a global pharmaceutical supply chain design model that considers the risk of export bans that are induced by supplier capacity disruptions and corresponding price increases. The model takes the company's perspective as a decision-maker looking to locate plants and distribute drugs globally. It is a two-stage stochastic program that includes uncertainty in capacity, ability-to-export, and demand. The model is solved by integrating the Sample Average Approximation and L-shaped methods. We present conditions related to when demand will be met and a case study of a generic oncology drug. We find that preparing for geopolitical strain may increase resilience and profits as well as reduce shortages in the short term. At baseline, expected global shortages are high (17.2\%) with disparities across country income levels (0.3\%, 0.8\%, 87.2\%, and 87.6\% for high, upper-middle, lower-middle, and low income countries, respectively). Pricing policies may improve drug access overall, back-shoring may slightly improve access for the country where it is implemented, and bilateral alliances may not be effective at improving access. \\
	\end{abstract}
			
	\noindent%

	\spacingset{1.5}

\section{Introduction} \label{s:intro}

Pharmaceutical companies operate globally to meet worldwide demand for drugs. Planning in a global context exposes supply chains to periods of geopolitical instability that have become frequent and prolonged \citep{Kalish2021}. Ongoing strains include the China-US trade war, Japan-South Korea trade dispute, COVID-19 pandemic, and the Russia-Ukraine war \citep{Kalish2021,Simon2022}. Under such events, governments may take actions to protect national security that affect global supply chain performance, including restricting exports \citep{Simon2022,WTO2020-ER}. 

Recently, several countries have banned or threatened to ban the export of drugs in short supply, e.g., India, France, Poland, Greece, Norway, Spain, and Bulgaria \citep{WTO2022b,Bochenek2018,EscoAguiar2021}. This means that when there are capacity disruptions, companies not only mitigate the lower capacity to produce, they may also need to consider \textit{induced} disruptions, where drugs that are produced may not be able to leave a country. It is not clear how companies could respond or mitigate these new dynamics.  

In this paper, we develop a new optimization approach to global drug supply chain planning. \textcolor{mycolor}{It considers} how a pharmaceutical company may design a global supply chain when there are risks of export bans concurrent with other disruption and operational risks. We further consider the impacts of policies from external entities, e.g., governments and non-profits, who seek to reduce worldwide drug shortages, a byproduct of the company's global supply chain decisions.

Broadly, the pharmaceutical supply chain is comprised of suppliers, manufacturers, distribution centers, healthcare systems, and patients \citep{NASEM2022}. Production and distribution are dispersed worldwide and exposed to several governments' regulations \citep{Marques2020}, risks of disruptions and uncertainties \citep{Hasani2016}, and the political and economic stability of stakeholders \citep{Snyder2016}. Given the rigidity and complexity of such regulatory frameworks and the large capital investment required in this industry \citep{Zhao-tutorial}, facility location decisions cannot be easily changed once executed. 

There are many pharmaceutical companies in general, but few produce specific drugs. This phenomenon, known as industrial concentration, implies that individual company decisions can have widespread effects on global availability, e.g., Teva Pharmaceutical's decision not to produce vincristine led to shortages \citep{Dyerl6086}. In addition, the production of Active Pharmaceutical Ingredients (APIs) and some drugs is geographically concentrated; e.g., India has 20\% of global exports of generic drugs and relies on China for most of its APIs \citep{McKinseyGlobalInstitute2020}. Industrial and geographic concentrations increase supply chain vulnerabilities to geopolitical, regulatory, economic, and climate country issues \citep{TheWhiteHouse2021}, which could lead to global shortages.

The geographic concentration of drug manufacturing in the context of export bans affects the global economy, increasing drug prices and exacerbating the risk of global drug shortage \citep{Casey2021}. Low and middle income countries with limited local production and low power to compete for drugs in short supply are vulnerable and have difficulties in drug access \citep{Tatar2022}. Drug shortages affect the quality and continuity of patients' treatments, which endangers the patient's health and increases the costs of healthcare \citep{Bochenek2018}. Companies with manufacturing plants in countries with these trade barriers may suffer economic losses due to the inability to sell their products in attractive global markets.

Locating manufacturing plants in domestic countries (known as on-shoring or back-shoring) and international alliances have been suggested as solutions to emerging geopolitical strains \citep{AAM2021,TheWhiteHouse2021}. Governments from many countries have called for on-shoring manufacturing activities, e.g., Japan, France, and the US \citep{Barbieri2020,TheWhiteHouse2021}. Countries with limited manufacturing capacities are focusing on strategies to strengthen their production capacities, e.g., the Latin America and Caribbean Region \citep{ECLAC2021}. However, on-shoring decisions are not always economically feasible for companies. Building international alliances between nations may mitigate the risk of export bans by the premise that alliances reduce uncertainty on partners' behaviors \citep{Varadarajan1995}, but effects on supply chain design decisions and the company's performance are unknown.

Identifying and including uncertainties and disruptions in global supply chain design models, even outside of export ban risks are challenging in terms of modeling and solution approaches \citep{Sabouhi2018,Blossey2022,Marques2020}. It requires a careful parameterization of the uncertainty sets and the study of how those uncertainty elements affect the structure of the mathematical model, leading to large, complex stochastic programs. The integration of export policies, pricing, and regulations in global pharmaceutical supply chain operations remains unexplored \citep{Goodarzian2020,Diaz2022}. No studies have focused on helping pharmaceutical companies design their supply chains in an environment with increasing risks of export restrictions, nor have they explored effects on drug availability across income levels.

We aim to study the pharmaceutical supply chain design problem in a global context. Under this approach, we include uncertainties and risks external to the supply chain network, e.g., geopolitical strain; internal to the company, e.g., facility quality issues; and external to the company but internal to the network, e.g., uncertainty in demand and availability of raw materials. Our goal is to produce insights into how emerging export ban risks impact pharmaceutical supply chain design decisions, economic performance, and drug shortages in a mix of markets \textcolor{mycolor}{(high, upper-middle, lower-middle, and low-income countries)}. We summarize the contributions of our study as follows:
\begin{itemize}[noitemsep,topsep=0pt]
    \item \textcolor{mycolor}{We introduce} the first pharmaceutical supply chain design model that incorporates export bans and export ban-induced price increases as a manifestation of geopolitical strain.
    \item We analyze the effects of external stakeholders' policies oriented to improve drug access on company economic performance, design decisions, and drug shortages globally and by income level classification of countries. These include exogenous bilateral alliances between countries, pricing, and back-shoring policies.
\end{itemize}

The remainder of this study is as follows. In Section \ref{s:lit-rev}, we review relevant literature. In Section \ref{s:model}, we present our model and its structural properties, and in Section \ref{s:methods}, we present the solution methods. Section \ref{s:case} evaluates the model in an oncology drug case study, presents several computational experiments, and discusses our results. In Section \ref{s:conclusion}, we present our conclusions and future directions.

\section{Literature review} \label{s:lit-rev}
Our research relates to the literature on supply chain design models under uncertainty, with a focus on facility location. In this section, we discuss supply chain design models generally (from any industry) and then concentrate on pharmaceutical supply chains. 

Global models consider a multi-national context for design decisions. A review is available in \citet{KchaouBoujelben2018}; it indicates that there are few global facility location models with uncertainty and that for-profit models focus on international economic features, e.g., transfer prices. Uncertainty is considered in demand, exchange rates, taxes, costs, and prices \citep{Goh2007,KchaouBoujelben2018}. In not-for-profit settings such as humanitarian logistics, decisions include opening warehouses around the world and pre-positioning relief items \citep{Duran2011,Jahre2016}. Both for-profit and not-for-profit settings leave out geopolitical stability issues that can emerge in globalized contexts, and for-profit models do not analyze disparities between countries in product access. Empirical research analyzes drivers for on-shoring (domestic country) and off-shoring (foreign country) manufacturing capacities \citep{Mohiuddin2019}. \citet{Hilletofth2019} suggest that instead, right-shoring (a mix between domestic and foreign locations) can improve upon both. Right-shoring decisions are made to support competitiveness and customer service and are based on factors including uncertainty in demand, supply, costs, disasters, and country risks \citep{McIvor2021}.

\textcolor{mycolor}{There is a wide literature on facility location, and} we review those with features and uncertainty most relevant to this work. Optimization models commonly consider uncertainty in demand and capacity \citep{Govindan2017,Tordecilla2021}. Facility capacity is predominantly modeled as binary availability, where a facility's ability to produce is either available or completely disrupted; \textcolor{mycolor}{a survey on models that include disruptions is available in \citet{Snyder2016}}. Recent studies include partial disruptions of capacity \citep{Aldrighetti2021, Maharjan2022}. For instance, \citet{Ghavamifar2018} consider partial \textcolor{mycolor}{effects} on the capacity of distribution centers, and \citet{Namdar2018} on suppliers' capacities. \textcolor{mycolor}{In all of these cases,} capacity is modeled as a \textcolor{mycolor}{single generic parameter that is uncorrelated with other modeled uncertainty and does not affect prices.}

\textcolor{mycolor}{Opportunity costs for unmet demand, also considered as penalty costs or shortage costs, are generally incorporated with a fixed, deterministic coefficient} \citep{Tordecilla2021}. One exception models stochastic penalty costs with a uniform distribution \citep{Namdar2018}. Shortage costs may also be deterministic but increased compared to a pre-disruption baseline. For example, in the disaster relief literature, shortage costs are based on an increase in the commodity's price (e.g., ten times the procurement price \citep{Rawls2010}). There is also research related to supply chain planning that includes economic conditions that may result in price changes, e.g., \citet{AZAD2019} propose a model for selecting optimal recovery strategies from supply disruptions that considers a dynamic pricing strategy as a mechanism to manage a price-sensitive demand. Their model assumes that the firm has market power.

While geopolitical and social risks, including regulatory issues, trading barriers, and strikes, are recognized as potential disruptions with serious consequences \citep{Kochan2018,Ambulkar2015}, they receive little consideration in the facility location literature \citep{Chatzoglou2018, Suryawanshi2022}. In one example, \citet{Nguyen2021} evaluate the vulnerability of a multi-echelon assembly supply chain network to different disruptive events, including labor strikes. New models and approaches are needed to parameterize the growing risks of export bans.

Research about alliances and facility location models is similarly limited. These partnerships may reduce uncertainty in the business environment and partners' behavior \citep{Das2018}, and empirical studies focus on determining the drivers of alliances \citep{Singh2018}. They are studied in graph theory \citep[e.g.,][]{Ouazine2018} and game theory \citep[e.g.,][]{LiXiaopeng2020}. Within facility location, \citet{Namdar2018} present an example of supplier-buyer alliances. They propose a model for risk mitigation of supply shortages, which integrates collaboration and visibility to improve suppliers' recovery capabilities and buyer's warning capability. \citet{Jahre2016} highlight the need for bilateral agreements with governments to facilitate disaster relief operations. Supply chain partnerships may support resilience, including managing risks \citep{Singh2019}.

Next, we focus on research related to pharmaceutical supply chains. \citet{NASEM2022} proposes a framework to improve resilience that includes a mix of mitigation, preparedness, response, and awareness interventions. They recommend international treaties and diversification of supply chains to protect patients from supply interruptions. Modeling studies mainly focus on optimizing operations and inventory management \citep{Franco2017,Blossey2021}. Those that include uncertainty generally consider uncertainty in demand \citep{Blossey2021}. Regulatory uncertainty, such as market authorization \citep{HANSEN2015}, is also considered.

Among supply chain design models, recent studies consider uncertainty and disruptions. \citet{Tucker2020} propose stochastic programs that include disruptions in multiple echelons and recovery time and study the effects of domestic policies to reduce drug shortages. \citet{LiJinfeng2023} propose a distributionally robust optimization model that considers uncertainty in demand and non-shortage costs. \citet{Tucker2022} develop closed-form expressions for supply chain reliability. Some studies focus on particular aspects of pharmaceuticals, such as perishability, e.g., \citet{Zandkarimkhani2020}; uncertainty is generally considered in demand and cost (fixed, inventory, and transportation).

\textcolor{mycolor}{There are very few studies that consider pharmaceutical planning in uncertain multi-national settings. \citet{Hasani2016} consider three international markets and propose a non-linear model with a hybrid parallel meta-heuristic to optimize a supply chain that is robust against correlated disruptions and uncertainty in procurement costs and demand. \citet{Goodarzian2020} also take a meta-heuristic approach to solve a multi-objective supply chain design model with multi-modal transportation systems.} The model allows investment in foreign manufacturing centers to satisfy local demand in hospitals and pharmacies. Uncertainty is considered in costs and capacities. \textcolor{mycolor}{\citet{Blossey2022} propose a multi-period and multi-product stochastic program with uncertainty in regulatory approval times and demand to optimize capacity expansion, product allocation, and selection of contract manufacturing organizations.} The numerical study is based on five regions. In a discrete-event simulation approach, \citet{Diaz2022} design a framework that assesses logistics drivers' effects on global pharmaceutical supply chains, considering uncertainty in demand, supply, and production. \textcolor{mycolor}{Deterministic global models have incorporated} economic features of international trade such as taxes \citep{Sousa2011}, import duties, and transfer prices \citep{Susarla2012}. \textcolor{mycolor}{There is not yet a stochastic model that considers the full global context, and none of the current approaches incorporate geopolitical strain.}

Recent studies integrate \textcolor{mycolor}{alliance} structures between the private players of the supply chains. \citet{Iacocca2022} develop an optimization model to know how cooperative drug-specific partnerships between mail-order and chain pharmacies can be used to reduce drug prices and increase profits. \citet{Akbarpour2020} propose a bi-objective model with supplier agreements and locate mobile pharmacies within a relief network design according to the cooperative coverage mechanism to respond to demand. \citet{LiZhao2023} propose a competition and cooperation model between pharmaceutical and third-party logistics enterprises for drug distribution. No studies consider alliances between governments as stakeholders of this supply chain.

From the perspective of mitigating drug shortages, there is research on inventory management in healthcare facilities (e.g., \citet{SAEDI2016}), as well as in the area of government interventions. Recent frameworks have been developed to support policy-makers in evaluating the impacts of government regulations. \citet{DONG-LI2022} evaluate the impacts of price and purchase regulations with a game theoretic model; \citet{Heese2023} determine the effect of interventions in the total industry equilibrium supply, and \citet{Kazaz2023} analyze the effects of social investor' interventions based on subsidies for encouraging manufacturers to build production and distribution capacity to serve low and middle income markets. Note that these modeling frameworks have been from third-party agents' perspectives, either governments or healthcare facilities, instead of the pharmaceutical company as the decision-maker. In addition, they do not consider a mix of markets, such as high, middle, and low income level markets, which hide dynamics that can emerge due to the global features of the pharmaceutical supply chain. Interventions or instruments used by social investors (governments or philanthropic institutions) oriented to improve drug access in low and middle income countries may also have effects in high income countries. 

Taken together, there is rich literature on supply chain and pharmaceutical planning. However, current frameworks are not sufficient to address key questions related to geopolitical strain. \textcolor{mycolor}{It is unknown how geopolitical risks may influence the supply chain design problem and shortages globally. In particular, the behavior of export bans is different from other types of disruptions; the risk of bans is activated by other types of supply strain, i.e., supplier capacity, and imposed bans increase prices, differentiated by domestic vs. exported supply.} Existing approaches in the supply chain disruption literature are not sufficient to capture these dynamics. This paper contributes to the literature with a modeling approach that considers global supply chains, \textcolor{mycolor}{geopolitical risks of} export bans and export ban-induced price increases, \textcolor{mycolor}{capacity} strains and disruptions, and alliances between governments \textcolor{mycolor}{that may mitigate} geopolitical strain.

\section{Model description} \label{s:model}
To consider how geopolitical strain may affect a company's plans and drug shortages, we introduce a new global pharmaceutical supply chain model for a single drug. The model takes the perspective of the pharmaceutical company as the decision-maker. It is a two-stage stochastic program that aims to locate manufacturing plants and decide material flows in the presence of uncertainty, \textcolor{mycolor}{maximizing the company's expected profits}. The model is a three-echelon supply chain formed by suppliers $(I)$, plants $(J)$, and demand countries $(K)$. \textcolor{mycolor}{Each country $k$, has a set of allied $(F_{k})$ and a set of non-allied $(F^{\prime}_{k})$ countries.}

The uncertainty elements are geopolitical strain (export bans and export ban-induced price increases), facility operational uncertainty (strains and disruptions in the capacity of suppliers and plants), and demand. Exogenous, \textcolor{mycolor}{bilateral alliances between $k$ and other nations ($F_{k}$)} affect export ban uncertainty to the allied countries. To consider production of the drug by other stakeholders, we also model exogenous drug exports from each country. The notation is presented in Table \ref{Table:Notation}.

\begin{figure}[ht]
    \centering
    \includegraphics[width=15cm]{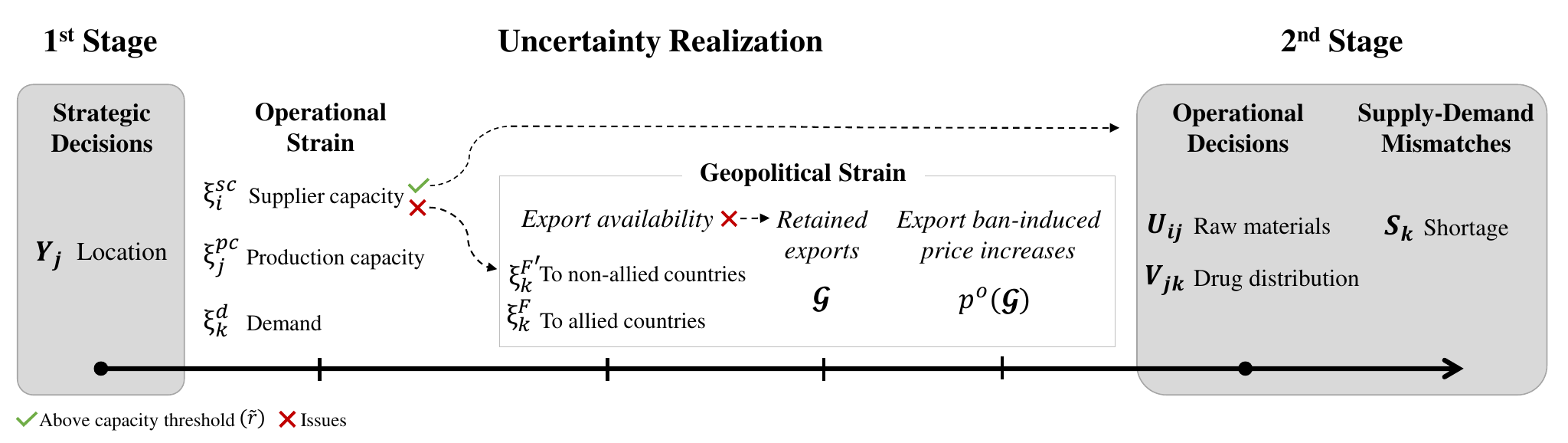}
    \caption{\textcolor{mycolor}{Sequence of events.}}
    \label{fig:sequence-of-events.eps}
    \vspace{-10pt}
\end{figure}

The sequence of events can be seen in Figure \ref{fig:sequence-of-events.eps}. In the first stage, manufacturing plants are located. In the second stage, after uncertainty is realized, the drug is produced, distributed, and there may be shortages. Between the stages, the suppliers' capacities $(\xi_i^{sc})$, production capacities $(\xi_j^{pc})$, and demand $(\xi_k^{d})$ are first realized. Then, based on the \textcolor{mycolor}{weighted} average global availability of raw material capacity, export bans are implemented considering the alliance status, \textcolor{mycolor}{$(\xi_k^{F^\prime})$} for non-allies, and \textcolor{mycolor}{$(\xi_k^{F})$} for allies. Finally, an increase in prices \textcolor{mycolor}{$(p^o)$} may occur based on the retained exports generated by the export bans $(\mathcal{G})$. 

\begin{table}[ht]  
\scriptsize
\renewcommand{\arraystretch}{0.85}
\setlength{\tabcolsep}{10pt}
  \caption{Notation.}
\begin{tabular}{ll}
\hline
\multicolumn{2}{l}{\textbf{Sets}}  \\
$I$ & set of supplier countries.\\ 
$J$ & set of potential countries to locate a manufacturing plant.\\
$K$ & set of demand countries. $I\subseteq K$ and $J\subseteq K$.\\
\textcolor{mycolor}{$F_{k}$} & \textcolor{mycolor}{set of allied countries of country $k \in K$. Note $k \not \in F_{k}$.}\\
\textcolor{mycolor}{$F_{k}^{\prime}$} & \textcolor{mycolor}{set of countries that are non-allies of country $k \in K$. Note $F_{k}^{\prime} = K\setminus (F_{k}\cup\{k\})$.}\\
$\Omega$ & set of scenarios (discrete).\\
\multicolumn{2}{l}{\textbf{Parameters}} \\
\textcolor{mycolor}{$\mathcal{P}^\omega$} & probability of scenario $\omega \in \Omega$. \\  
$c_{i}^{rm}$ & unit cost of raw material purchased from supplier $i \in I$. \\ 
$c_{j}^{pr}$ & unit production cost of manufacturing plant $j \in J$. \\
$c_{j}^{fi}$ & annual cost for opening and operating manufacturing plant $j \in J$. \\
$c_{ij}^{t1}$ & unit transportation cost from supplier $i \in I$ to manufacturing plant $j \in J$. \\
$c_{jk}^{t2}$ & unit transportation cost from manufacturing plant $j \in J$ to country $k \in K$. \\
$\textcolor{mycolor}{p}_k^s$ & baseline unit drug price in country $k \in K$.\\
$q_{i}^{sc}$ & capacity of supplier $i \in I$. \\
$q_{j}^{pc}$ & production capacity of manufacturing plant $j \in J$.\\
\textcolor{mycolor}{$e_{k}^{F}$} & \textcolor{mycolor}{exogenous drug exports from country $k \in K $ to allied countries} \\
\textcolor{mycolor}{$e_k^{F^\prime}$} & \textcolor{mycolor}{exogenous drug exports from country $k \in K$ to non-allies.}\\
$\beta$ & coefficient for export ban-induced price increases. \\
$\rho_k$ & probability of country $k \in K$ allowing exports \textcolor{mycolor}{to allies}.\\
\textcolor{mycolor}{$\rho^{\prime}_k$} & probability of country $k \in$ \textcolor{mycolor}{$K$ allowing exports to non-allies}.\\
\multicolumn{2}{l}{\textbf{Random Parameters}} \\
$\xi_i^{sc,\omega}$ & proportion of capacity of supplier $i \in I$ that is available in scenario $\omega \in \Omega$.\\
$\xi_j^{pc,\omega}$ & proportion of capacity of manufacturing plant $j \in J$ that is available in scenario $\omega \in \Omega$.\\
$\xi_k^{d,\omega}$ & demand $k \in K$ in scenario $\omega \in \Omega$.\\
\textcolor{mycolor}{$\xi_k^{F,\omega}$} & binary variable that equals 1 if nation $k \in K$ allows exports \textcolor{mycolor}{to allies} and 0 otherwise.\\
\textcolor{mycolor}{$\xi_k^{F^{\prime}\mspace{-7mu},\omega}$} & binary variable that equals 1 if nation $k \in$ \textcolor{mycolor}{$K$} allows exports \textcolor{mycolor}{to non-allies} and 0 otherwise.\\
$\mathcal{G}^\omega$ & global retained exports in scenario $\omega \in \Omega.$\\
$\textcolor{mycolor}{p}^{o,\omega}$ & marginal unit drug price increase due to export bans in scenario $\omega \in \Omega$.\\
\multicolumn{2}{l}{\textbf{Decision Variables}} \\
$Y_{j}$ & binary variable that equals 1 if location $j \in J$ is selected and 0 otherwise.\\
$U_{ij}^{\omega}$ & raw material purchased by manufacturing plant $j \in J$ from supplier $i \in I$ in scenario $\omega \in \Omega$. \\
$V_{jk}^{\omega}$ & drugs produced by manufacturing plant $j \in J$ sent to country $k \in K$ in scenario $\omega \in \Omega$. \\
$S_{k}^{\omega}$ & drug shortage in country $k \in K$ in scenario $\omega \in \Omega$.\\
\hline
\end{tabular}\label{Table:Notation}
\vspace{-10pt}
\end{table}

More specifically, the production capacity of plants $(\xi_j^{pc})$ is defined by two types of facility interruptions: strain $(\xi_j^{ps})$ and disruptions $(\xi_j^{pd})$. Strain refers to partial capacity interruptions; \textcolor{mycolor}{we model it as a discrete random variable with probability mass function $p(n)$, where $n$ represents capacity levels}. Disruption refers to the complete unavailability of the facility; we model it as $\xi_j^{pd}$ $\sim {\rm Bernoulli}(\rho_j^{pd})$ for all $j \in J$, where $\rho_j^{pd}$ \textcolor{mycolor}{is} the probability of having facility $j$ available. The combined effect of those two interruptions is $\xi_j^{pc} = \xi_j^{ps} \xi_j^{pd}$. The raw material available at each supplier ($\xi_i^{sc}$) is also defined by the two types of interruptions above (strains $\xi_i^{ss}$ and disruptions $\xi_i^{sd}$), modeled like the capacity of plants. 

The model considers the risk of export bans in every country. In practice, export bans could be implemented for several reasons; in this model, the initiating condition is based on the \textcolor{mycolor}{weighted} average global capacity for raw material production \textcolor{mycolor}{$\big(\sum_{i \in I}q_{i}^{sc}\xi_{i}^{sc}/\sum_{i \in I}q_{i}^{sc}\big)$}. If \textcolor{mycolor}{$\big(\sum_{i \in I}q_{i}^{sc}\xi_{i}^{sc}/\sum_{i \in I}q_{i}^{sc}\big)$} is less than a threshold $\tilde{r}$, there is a leading indicator of a risk of shortage, \textcolor{mycolor}{and} then there exists a risk that governments implement export bans; otherwise, there is no geopolitical strain. 

The model considers exogenous bilateral alliances between \textcolor{mycolor}{country $k \in K$ and other countries $F_{k}$}. Members of this group present less risk of being affected by their partner country's export bans. This feature requires two random parameters related to export bans, \textcolor{mycolor}{$\xi_k^{F^\prime}$} for non-allies and \textcolor{mycolor}{$\xi_k^{F}$} for allies. We sample \textcolor{mycolor}{$\xi_k^{F^\prime}$} $\sim {\rm Bernoulli}$ (\textcolor{mycolor}{$\rho_k^\prime$}) for all $k \in K$ and apply for non-allies of $k$ \textcolor{mycolor}{($F_k^{\prime}$)}, and if \textcolor{mycolor}{$\xi_k^{F^\prime}$}$=0$, then sample \textcolor{mycolor}{$\xi_k^{F}$} $\sim {\rm Bernoulli}(\textcolor{mycolor}{{\rho_k}})$ \textcolor{mycolor}{and apply to the set of allies ($F_k$)}. The terms \textcolor{mycolor}{$\rho_k^\prime$ and $\rho_k$} represent the probability of allowing exports. Note that having alliances may mitigate export ban uncertainty. To illustrate, if a country $k$ allows exports to its non-allied countries (\textcolor{mycolor}{$\xi_k^{F^{\prime}\mspace{-7mu},\omega}$} = 1), then it allows exports to its allies too (\textcolor{mycolor}{$\xi_k^{F,\omega}$} = 1). If a country $k$ bans export to its non-allied countries (\textcolor{mycolor}{$\xi_k^{F^{\prime}\mspace{-7mu},\omega}$} = 0), $k$ may not restrict exports to its allies (\textcolor{mycolor}{$\xi_k^{F,\omega}$} can be 0 or 1). The model assumes that the joint probability distribution of the random parameters $\xi$ $:=$ $[ \xi^{sc}, \xi^{pc}, \xi^{d}, \textcolor{mycolor}{\xi^{F^\prime}, \xi^{F}}]^\top$ is known.

By reducing the global drug supply, export bans may increase drug prices. We model this decrease as the global retained exports, $\mathcal{G}^\omega$ (Equation \eqref{eq:calG}). This value represents the total exogenous drug supply (\textcolor{mycolor}{$e_k^{F^\prime}, e_k^{F}$}) that is prevented from leaving countries due to export bans in scenario $\omega \in \Omega$. \textcolor{mycolor}{The first and second terms refer to global exports to non-allies and allies, respectively.}
\begin{equation}\label{eq:calG}
    \textcolor{mycolor}{\mathcal{G}^{\omega} = \sum_{k \in K} e_k^{F^\prime}(1-\xi_k^{F^{\prime}\mspace{-7mu},\omega}) + e_k^{F}(1-\xi_k^{F,\omega}) } 
\end{equation}

\textcolor{mycolor}{Drug prices are modeled with two elements: the baseline price of the drug ($p_k^s$) and the stochastic marginal price increase due to export bans ($p^{o,\omega}(\mathcal{G}^\omega)$). Note that the marginal price increase is a function of retained exports during bans, $\mathcal{G}^\omega$. For simplicity, we apply a linear relationship, i.e., $p^{o,\omega}= \beta \mathcal{G}^\omega$, though other functions could be incorporated without loss of generality. If there are no export bans, the drug price in each country $k\in K$ is set to its baseline value ($p_k^s$). If there are export bans in any country, most prices are the baseline value plus the marginal price increase ($p_k^s+p^{o,\omega}$). The only exception is for domestic sales in countries that impose bans; these prices remain at baseline because these sales are unaffected by the global market \citep{WTO2020-ER}.}

\begin{figure}[h]
    \centering
    \includegraphics[width=9cm]{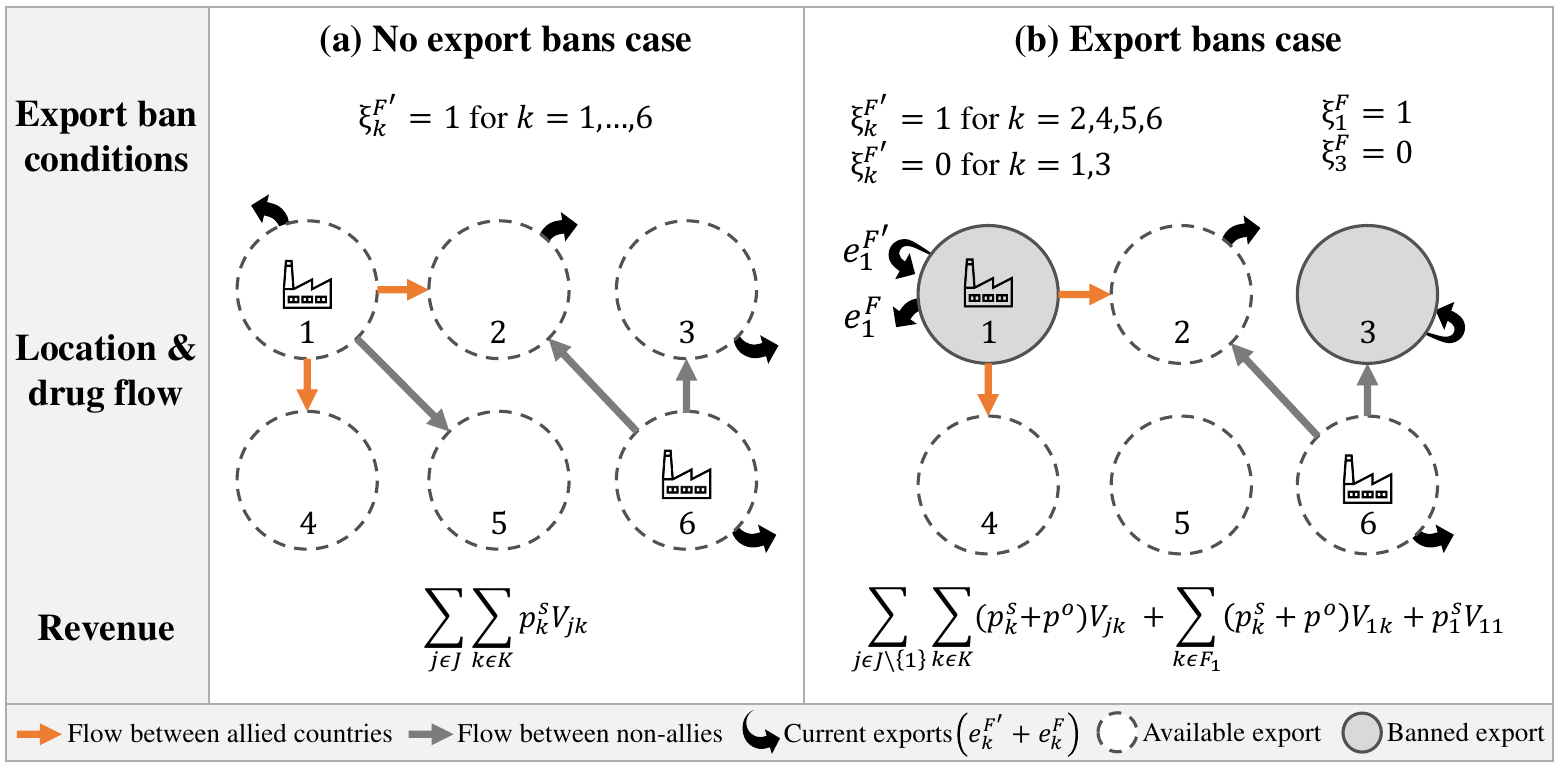}
    \caption{\textcolor{mycolor}{Revenue composition.}}
    \label{fig:opportunity-cost.eps}
    \vspace{-10pt}
\end{figure}

\textcolor{mycolor}{Figure \ref{fig:opportunity-cost.eps} illustrates the price mechanics. Without bans (panel a), the prices for sales to each country are at the baseline prices ($p_k^s$). In panel (b), there are bans; country $k=1$ bans exports to non-allies ($\xi_1^{F'} =0$) and allows exports to allies ($\xi_1^F=1$), and country $k=3$ bans exports to all countries. The total retained exports are $e_{1}^{F^\prime}+e_{3}^{F^\prime}+e^{F}_{3}$, and the marginal price increase is $p^o = \beta ( e_{1}^{F^\prime}+e_{3}^{F^\prime}+e^{F}_{3})$. The price of the company's sales for domestic production under an export ban ($k=1$) are set to baseline $p_1^s$. All other prices are influenced by the drop in global supply and are set to the baseline plus the marginal price increase.}

\subsection{Model formulation} \label{s:formulation}
The formulation of the stochastic program is as follows: 
\allowdisplaybreaks
\begin{subequations}\label{2SSP}
    \small
        \begin{align}
            z= & \textcolor{mycolor}{\mspace{5mu} \max \limits_{Y} \left \lbrace -\sum_{j \in J} c_j^{fi}Y_j + \sum_{\omega \in \Omega} \mathcal{P}^\omega Q^{\omega}(Y) \right \rbrace }\label{OF1}\\
             \text{s.t.  } & \sum_{j \in J}Y_j \geq 1 \label{C1.1}\\
             & Y_j \in \{0,1\} \quad \forall j \in J \label{C1.2}
        \end{align}
        \vspace{-10mm}
        \begin{align}
            \textcolor{mycolor}{Q^{\omega}(Y):= }& \textcolor{mycolor}{\mspace{5mu} \max \sum_{j \in J} \sum_{k \in K} (p_{k}^{s}-c_j^{pr}-c_{jk}^{t2})V_{jk}^{\omega} + p^{o,\omega} \sum_{j \in J} \left( \sum_{k \neq j \in K} V_{jk}^{\omega} + \xi_{j}^{F^{\prime}\mspace{-7mu},\omega}V_{jj}^{\omega} \right) - \sum_{i \in I}\sum_{j \in J}(c_i^{rm}+c_{ij}^{t1})U_{ij}^\omega} \label{OF2}\\
             \text{s.t.  } & \sum_{j \in J}U_{ij}^\omega \leq q_{i}^{sc}\xi_i^{sc,\omega} \quad \forall i \in I \label{C2.1}\\
             & U_{ij}^\omega \leq q_{i}^{sc}\xi_i^{sc,\omega}\textcolor{mycolor}{\xi_i^{F^{\prime}\mspace{-7mu},\omega}}Y_j \quad \forall i \in I, \textcolor{mycolor}{j \in (J \cap F_{i}^{\prime})} \label{C2.2}\\
             & U_{ij}^\omega \leq q_{i}^{sc}\xi_i^{sc,\omega}\textcolor{mycolor}{\xi_i^{F,\omega}}Y_j \quad \forall i \in I, \textcolor{mycolor}{j \in (J \cap F_{i}) } \label{C2.3}\\
             & \sum_{k \in K}V_{jk}^\omega \leq q_{j}^{pc}\xi_j^{pc,\omega}Y_j \quad \forall j \in J \label{C2.4}\\
             & V_{jk}^\omega \leq q_{j}^{pc}\xi_j^{pc,\omega}\textcolor{mycolor}{\xi_j^{F^{\prime}\mspace{-7mu},\omega}}Y_j \quad \forall j \in J, \textcolor{mycolor}{k \in F_{j}^{\prime}} \label{C2.5}\\
             & V_{jk}^\omega \leq q_{j}^{pc}\xi_j^{pc,\omega}\textcolor{mycolor}{\xi_j^{F,\omega}}Y_j \quad \forall j \in J,\textcolor{mycolor}{k \in F_{j}} \label{C2.6}\\
             & \textcolor{mycolor}{\sum_{j \in J}V_{jk}^\omega + S_k^\omega = \max \big \{\xi_k^{d,\omega} - (1-\xi_k^{F^{\prime}\mspace{-7mu},\omega})e_{k}^{F^{\prime}} - (1-\xi_k^{F,\omega})e_k^{F},0 \big \} \quad \forall k \in K} \label{C2.7}\\
             & \sum_{i \in I}U_{ij}^{\omega} = \sum_{k \in K}V_{jk}^\omega \quad \forall j \in J \label{C2.8}\\
             & \textcolor{mycolor}{U_{ij}^{\omega}, V_{jk}^{\omega}, S_{k}^{\omega} \geq 0 \quad \forall i \in I, j \in J, k \in K} \label{C2.9}
        \end{align} 
\end{subequations}

The objective is to \textcolor{mycolor}{maximize expected profits. In the first stage, the objective function \eqref{OF1} maximizes the negative fixed costs plus the total expected profits of the sale of drugs worldwide}. Constraint \eqref{C1.1} ensures that at least one plant is open, and constraints \eqref{C1.2} enforce the domain. In the second stage, for a finite set of scenarios $\Omega$, it \textcolor{mycolor}{maximizes the total profits \eqref{OF2}, considering revenue, costs of purchases of raw materials, production and transportation. To calculate the revenue, two types of prices are considered, $p_{k}^s$, the price of the drug in the country $k$, and $p^{o,\omega},$ the marginal increase in prices due to the retained exports.} Raw materials purchases are limited by the available supplier capacity \eqref{C2.1}. Flows of raw materials are limited by export bans; constraints \eqref{C2.2} limit flows between non-allied countries and \eqref{C2.3} between allies. Similarly, the production and distribution of final products are limited by the available manufacturing capacity \eqref{C2.4} and export bans, where constraints \eqref{C2.5} limit flows between non-allied countries, and \eqref{C2.6} between allies. The demand can be satisfied by the company's production or by the exports held in the country because of an export ban \textcolor{mycolor}{(constraints \eqref{C2.7}).} Constraints \textcolor{mycolor}{\eqref{C2.8}} correspond to the flow balance constraints. Finally, domain constraints for the second-stage variables are enforced \textcolor{mycolor}{\eqref{C2.9}}. 

\subsection{Model assumptions and justification} \label{s:m_assump}

This model is subject to several assumptions. We use a stochastic programming, expected value framework because the model represents the perspective of a for-profit pharmaceutical company \textcolor{mycolor}{that produces generic, injectable drugs. Margins to produce generic injectable drugs tend to be tight, and companies generally operate in a risk-neutral, rather than risk-averse context (personal communications with pharmaceutical manufacturers). Following these practices and the literature \citep{Tucker2020}, we apply a risk-neutral approach to maximize the modeled company's expected profit. A framework that considers robustness (e.g., a robust or distributionally-robust model) would be appropriate for extensions to a company that sells high-margin branded drugs or for a non-profit organization focused on drug access.} We consider the drug demand as price inelastic and uncorrelated because our work focuses on essential drugs that do not have substitutes \citep{NASEM2022,Zhao-tutorial}, and drug demand is mainly driven by the drug's efficacy and the number of patients \citep{Blossey2022}. The implementation of export bans is independent between countries but conditional on the risk of supplier capacity disruption; this is because we focus on the direct (first-order) effects of export bans instead of the retaliations that governments can take against each other due to trade barriers. Changes in prices produced by a mismatch between supply and demand other than those caused by export bans are not considered, given that we are modeling marginal capacities and demand, not the entire market. We assume that capacity strains and disruptions are independent, \textcolor{mycolor}{consistent with available literature \citep{Tucker2020}. We also assume that supplier and manufacturing disruptions last over a year; disruptions last 1.2 years on average at suppliers and 0.8 years on average at manufacturers \citep{Tucker2020}}. We consider static location decisions because of the inherent rigidity of pharmaceutical supply chains due to regulation \citep{Zhao-tutorial}. Since our focus is on strategic decisions taken by a single company, we exclude granular operational decisions (e.g., inventory and drug distribution within each market) and competition; other manufacturers are approximated via the exogenous export parameters. Lastly, the model focuses on the first-order effects of geopolitical strain by applying a single period in the second stage. Pharmaceutical export bans last over a year on average \citep{WTO2022a,WTO2022b}, and the fixed costs in the first stage are amortized in the case study to represent annual costs. The model then represents the expected annual \textcolor{mycolor}{profits} in the presence of geopolitical risks.

\subsection{Structural properties} \label{s:properties}
\textcolor{mycolor}{This section presents theoretical properties of the model, which help to understand its behavior and when demand may be met. These properties allow us to study the results of second-stage decision variables without solving problem \eqref{2SSP}.}

\textcolor{mycolor}{Theorem \ref{Thm:market} focuses on a country $k^\prime \in K$ imposing an export ban and states when retained exports will satisfy the country's demand, i.e., forcing $V_{j{k^\prime}}^{\omega}=0$ for all $j \in J$. This theorem reflects the company when there is an oversupply of the drug due to retained exports. In our model, the oversupply decreases the company-experienced demand due to an increase in other supply, i.e., retained exports and forces their sales to zero.}

\textcolor{mycolor}{Theorem \ref{Thm:cond_demand} states the necessary conditions for the company to satisfy demand in country $k^\prime \in K$ from manufacturing plant location $j^\prime \in J$ ($V_{{j^\prime}{k^\prime}}^{\omega} > 0$). In general, if there is production capacity in a manufacturing plant $j^\prime$, two things must occur for that plant to meet demand in country $k^\prime$. It must be economically feasible, and international trading must be allowed for flows between key countries. If at least one condition is not satisfied, then demand in that country will not met. This theorem provides a way to check whether a country may be at risk for not having its demand met. Parts \textit{(I)-(IX)} present the conditions for every combination of allied and non-allied categories of suppliers, plants and demand countries. To illustrate, consider part \textit{(I)}, where country $k^\prime$ is an ally of country $j^\prime$ (which is the country with the plant), and there exists an allied supplier $i$ (an ally of $j^\prime$) with capacity available. The plant $j^\prime$ will not meet demand in $k^\prime$ ($V_{j^{\prime}k^{\prime}}^{\omega}=0$) if at least one of the necessary conditions \textit{(i)}, \textit{(iii)a} and \textit{(iii)b} is not satisfied.} 

\textcolor{mycolor}{Finally, Theorem \ref{Thm:suf_demand} presents the sufficient condition for the company to meet demand, assuming that Theorem \ref{Thm:cond_demand} is satisfied. This theorem states that the company will meet demand in countries with higher marginal profits.}

\begin{theorem}[Sufficient conditions for $V_{j{k^\prime}}^{\omega}=0$ $\forall j \in J$]\label{Thm:market}
    \textcolor{mycolor}{For a given scenario $\omega \in \Omega$ and $k^\prime \in K$, if at least one condition (i)-(iii) is met, (i) \underline{ban for non-allies} $\xi_{k^\prime}^{F^{\prime}\mspace{-7mu},\omega}=0$ and $e_{k^\prime}^{F^{\prime}} \geq \xi_{k^\prime}^{d,\omega}$, or; (ii) \underline{ban for allies} $\xi_{k^\prime}^{F,\omega}=0$ and $e_{k^\prime}^{F}\geq \xi_{k^\prime}^{d,\omega}$, or; (iii) \underline{full ban} $\xi_{k^\prime}^{F^{\prime}\mspace{-7mu},\omega}=0$, $\xi_{k^\prime}^{F,\omega}=0$, and $(e_{k^\prime}^{F^{\prime}}+e_{k^\prime}^{F}) \geq \xi_{k^\prime}^{d,\omega}$. Then $S_{k^\prime}^{\omega}=0$ and $V_{j{k^\prime}}^{\omega}=0$ $\forall j \in J$.} 
\end{theorem}

\begin{proof}[\textbf{Proof}]
    \textcolor{mycolor}{Considering the theorem conditions and rewriting constraint (\ref{C2.7}), for $ k^{\prime} \in K,$ $\sum_{j \in J}V_{j{k^\prime}}^\omega + S_{k^\prime}^\omega  = 0$. By domain constraints (\ref{C2.9}), the only feasible solution is $\sum_{j \in J}V_{j{k^\prime}}^\omega = 0 \Rightarrow V_{j{k^\prime}}^{\omega}=0$ $\forall j \in J$, and $S_{k^\prime}^\omega=0$.}
\end{proof}

\begin{theorem}[Necessary conditions for $V_{{j^\prime}{k^\prime}}^{\omega} > 0$]\label{Thm:cond_demand}
Consider a scenario $\omega \in \Omega$, and conditions (i)-(iii), where:
\textcolor{mycolor}{
\begin{enumerate}[label=(\roman*),itemsep=0pt, topsep=0pt]
   \item Economic feasibility: $p_{k^\prime}^s+p^{o,\omega} > c_{i^\prime}^{rm}+c_{{i^\prime}{j^\prime}}^{t1}+c_{{j^\prime}}^{pr}+c_{{j^\prime}{k^\prime}}^{t2}$.
   \item International trade availability for non-allies: a. $\xi_{j^\prime}^{F^{\prime}\mspace{-7mu},\omega}=1$ and b. $\xi_{i^\prime}^{F^{\prime}\mspace{-7mu},\omega}=1$.
   \item International trade availability for allies: a. $\xi_{j^\prime}^{F,\omega}=1$ and b. $\xi_{i^\prime}^{F,\omega}=1$.
\end{enumerate} }

\vspace{5pt}
\noindent \textcolor{mycolor}{If $\exists j^\prime \in J$ such that $Y_{j^\prime}=1$ and $q_{j^\prime}^{pc}\xi_{j^\prime}^{pc,\omega}>0$, the necessary conditions for $V_{j^{\prime}k^{\prime}}^{\omega} > 0$ are:
\begin{enumerate}[label=(\Roman*),itemsep=0pt, topsep=0pt]
    \item For $k^\prime \in F_{j^\prime}$ and $i^\prime \in \{I \cap F_{j^\prime}|q_{i^\prime}^{sc}\xi_{i^\prime}^{sc,\omega}>0\}$, then (i), (iii)a, and (iii)b. 
    \item For $k^\prime \in F_{j^\prime}$ and $i^\prime \in \{I \cap F_{j^\prime}^{\prime}|q_{i^\prime}^{sc}\xi_{i^\prime}^{sc,\omega}>0\}$, then (i), (iii)a, and (ii)b.
    \item For $k^\prime \in F_{j^\prime}$ and $i^\prime=j^\prime$ with $i^\prime \in I$ and $q_{i^\prime}^{sc}\xi_{i^\prime}^{sc,\omega}>0$, then (i) and (iii)a.
    \item For $k^\prime \in F_{j^\prime}^{\prime}$ and $i^\prime \in \{I \cap F_{j^\prime}|q_{i^\prime}^{sc}\xi_{i^\prime}^{sc,\omega}>0\}$, then (i), (ii)a, (iii)b.
    \item For $k^\prime \in F_{j^\prime}^{\prime}$ and $i^\prime \in \{I \cap F_{j^\prime}^{\prime}|q_{i^\prime}^{sc}\xi_{i^\prime}^{sc,\omega}>0\}$, then (i), (ii)a, and (ii)b.
    \item For $k^\prime \in F_{j^\prime}^{\prime}$ and $i^\prime=j^\prime$ with $i^\prime \in I$ and $q_{i^\prime}^{sc}\xi_{i^\prime}^{sc,\omega}>0$, then (i) and (ii)a.
    \item For $k^\prime=j^\prime$ and $i^\prime \in \{I \cap F_{j^\prime}|q_{i^\prime}^{sc}\xi_{i^\prime}^{sc,\omega}>0\}$, then (i), (ii)a, and (iii)b.
    \item For $k^\prime=j^\prime$ and $i^\prime \in \{I \cap F_{j^\prime}^{\prime}|q_{i^\prime}^{sc}\xi_{i^\prime}^{sc,\omega}>0\}$, then (i), (ii)a, and (ii)b.
    \item For $k^\prime=j^\prime$ and $i^\prime=j^\prime$ with $i^\prime \in I$ and $q_{i^\prime}^{sc}\xi_{i^\prime}^{sc,\omega}>0$, then (i) and (ii)a.
\end{enumerate}
}
\noindent
\textbf{Note:} \textcolor{mycolor}{Under the no existence of export bans, condition (i) is harder to meet $(p_{k^\prime}^s > c_{i^\prime}^{rm}+c_{i{^\prime}j{^\prime}}^{t1}+c_{j^\prime}^{pr}+c_{j{^\prime}k{^\prime}}^{t2})$.}
\end{theorem}

\begin{proof}[\textbf{Proof}]
    \textcolor{mycolor}{Considering the premises presented for parts \textit{(I)-(IX)}, condition \textit{(i)} comes from demand constraints \eqref{C2.7}, domain constraints \eqref{C2.9}, flow balance constraints \eqref{C2.8} and the sense of the objective function \eqref{OF2}; conditions \textit{(ii)a} and \textit{(iii)a} comes from constraints \eqref{C2.5} and \eqref{C2.6}, respectively; and conditions \textit{(ii)b} and \textit{(iii)b} from constraints \eqref{C2.2} and \eqref{C2.3}, respectively. Note that flows with the same origin and destination, i.e., $V_{{j^\prime}{k^\prime}}^{\omega}$ for ${k^\prime}={j^\prime}$ and $U_{{i^\prime}{j^\prime}}^{\omega}$ for ${i^\prime}={j^\prime}$ are not included in \eqref{C2.5}-\eqref{C2.6} and \eqref{C2.2}-\eqref{C2.3}, respectively. Therefore, international trade availability conditions do not apply to those cases. However, for parts \textit{(VII)} to \textit{(IX), (ii)a} is required to keep condition \textit{(i)} valid due to objective function \eqref{OF2}.}
\end{proof}

\begin{corollary}\label{col-1-thm2}\textcolor{mycolor}{
\noindent For $k^\prime=j^\prime$, if condition (ii)a is not met, the necessary economic feasibility condition (i) is replaced by $p_{k^\prime}^s > c_{i^\prime}^{rm}+c_{i{^\prime}j{^\prime}}^{t1}+c_{j^\prime}^{pr}+c_{j{^\prime}k{^\prime}}^{t2}$ in parts (VII) to (IX).}
\end{corollary}

\begin{proof}[\textbf{Proof}]
    \textcolor{mycolor}{Since the objective function \eqref{OF2}, the final product flows between the same origin and destination ($V_{j^{\prime}j^{\prime}}^{\omega}$) are not subject to export ban-induced price increases $p^{o,\omega}$.}
\end{proof}

\begin{theorem}\label{Thm:suf_demand}
\textcolor{mycolor}{Consider a scenario $\omega \in \Omega$ and $Y_{j^\prime}=1$ for $j^\prime \in J$ with one unit of production capacity available. Suppose that for $k=k^\prime$ and $k=k^{\prime\prime}$, $\xi_{k}^{F^{\prime}\mspace{-7mu},\omega}=1$ and there exists $i^\prime \in \{I|q_{i^\prime}^{sc}\xi_{i^\prime}^{sc,\omega}>0\}$ such that Theorem \ref{Thm:cond_demand} is satisfied. Let $\delta_k^\omega$ be the unit profit of meeting demand in country $k$, i.e., $\delta_k^\omega:= (p_k^s+p^{o,\omega})-(c_{i^\prime}^{rm}+c_{{i^\prime}{j^\prime}}^{t1}+c_{j^\prime}^{pr}+c_{{j^\prime}k}^{t2})$. If $\xi_{k^\prime}^{d,\omega}>0$, $\xi_{k^{\prime\prime}}^{d,\omega}>0$ and $\delta_{k^\prime}^\omega > \delta_{k^{\prime\prime}}^\omega$, then $V_{j^\prime{k^\prime}}^\omega=1$ and $V_{j^\prime{k^{\prime\prime}}}^\omega=0$.}
\end{theorem}

\begin{proof}[\textbf{Proof}]
    \textcolor{mycolor}{Since Theorem \ref{Thm:cond_demand} is satisfied, we know that $V_{{j^\prime}{k^\prime}}^\omega$ and $V_{{j^\prime}{k^{\prime\prime}}}^\omega$ can be greater than zero, and that $\delta_{k^\prime}^\omega > \delta_{k^{\prime\prime}}^\omega > 0$. By the sense of the objective function \eqref{OF2} at optimality $V_{{j^\prime}{k^\prime}}^\omega=1$ and $V_{{j^\prime}{k^{\prime\prime}}}^\omega=0$.}
\end{proof}

\section{Solution approach} \label{s:methods}
To solve the model, we integrate the Sample Average Approximation (SAA) and L-shaped methods (Algorithm \ref{alg:SAA-BD}). SAA builds stochastic programs with scenario sets of a manageable size, and the L-shaped method solves them. \textcolor{mycolor}{This approach is well-suited for solving two-stage stochastic programs with continuous recourse \citep{Govindan2017}} and allows problems with large sample spaces to be solved with low optimality gaps in a reasonable time \citep{SANTOSO200596}. \textcolor{mycolor}{The latter is particularly relevant for the uncertainty in the considered problem; export ban risks, disruptions, operational strain, and demand are considered worldwide.}

\begin{algorithm}[ht] 
\LinesNumberedHidden
\spacingset{1.2}
\footnotesize
  \renewcommand{\arraystretch}{0.4}
    \caption{Integrated SAA and L-shaped decomposition methods}
      \label{alg:SAA-BD}
      \SetAlgoNoLine
      \DontPrintSemicolon
\KwIn{Calibrated numbers of $M$ replications, $N$ optimization scenarios and $N^\prime$ evaluation scenarios.
}
\KwOut{\textcolor{mycolor}{$L$ lower and $U$ upper bound estimators}, $\varepsilon$ optimality gap and $\hat{Y}^*$ $\varepsilon$-optimal solution.} 

\textbf{1.} In parallel \lFor{$m=1,...,M,$} {

$\quad  $ \textbf{1.1 } Generate $N$ i.i.d. scenarios of $\xi$.

$\quad  $ \textbf{1.2 } Solve the SAA problem \eqref{SAA_p} with the L-shaped method (Algorithm \ref{alg:BD}) and get $z_N^m$ and $\hat{Y}^m$.
\begin{align}\label{SAA_p}
    z_N^m= \textcolor{mycolor}{\max \limits_{Y}} \left \{ \textcolor{mycolor}{-\sum_{j \in J} c_j^{fi}Y_j^m + \frac{1}{N}\sum_{\omega = 1}^{N}Q^{\omega}(Y^{m})} 
    \right\} \\
    \text{s.t.} \quad \eqref{C1.1}-\eqref{C1.2}, \eqref{C2.1}-\textcolor{mycolor}{\eqref{C2.9}} \nonumber
\end{align}

$\quad  $ \textbf{1.3 } Generate $N^\prime$ i.i.d. scenarios of $\xi$.

$\quad  $ \textbf{1.4 } Use $\hat{Y}^m$ to evaluate the objective function $z_{N^\prime}^m(\hat{Y}^m)$, \textcolor{mycolor}{given by $-\sum_{j \in J} c_j^{fi}\hat{Y}_j^m + \frac{1}{N^\prime}\sum_{\omega = 1}^{N^{\prime}}Q^{\omega}(\hat{Y}^{m})$.}
}
\textbf{2.} Select $\hat{Y}^* \in \textcolor{mycolor}{\operatorname*{argmax}\{z^m_{N^\prime}(\hat{Y}):\hat{Y} \in \{\hat{Y}^1,...,\hat{Y}^M\}\}}$.

\textcolor{mycolor}{\textbf{3.} Calculate $U \coloneqq \bar{z}_{N,M}+t_{\alpha,\nu}\hat{\sigma}_{\bar{z}_{N,M}}$ upper bound of problem \eqref{2SSP}, where $\bar{z}_{N,M} \coloneqq \frac{1}{M} \sum_{m=1}^{M}z_N^m$, $\hat{\sigma}_{\bar{z}_{N,M}}^2 \coloneqq \frac{1}{(M-1)M} \sum_{m=1}^{M}(z_N^m-\bar{z}_{N,M})^2$, $\nu=M-1$, and $t_{\alpha,\nu}$ is the $\alpha$-critical value of the t-distribution with $\nu$ degrees of freedom.}

\textcolor{mycolor}{\textbf{4.} Calculate $L \coloneqq z_{N^\prime}-z_{\alpha}\hat{\sigma}_{z_{N^\prime}}$ lower bound of problem \eqref{2SSP}, where $z_{N^\prime}= -\sum_{j \in J} c_j^{fi}\hat{Y}^*_j + \frac{1}{N^{\prime}}\sum_{\omega = 1}^{N^{\prime}}Q^{\omega}(\hat{Y}^*)$ using a new sample of $N^\prime$ scenarios, and $\hat{\sigma}_{z_{N^\prime}}^2 \coloneqq \frac{1}{(N^\prime-1)N^\prime} \sum_{\omega=1}^{N^\prime}(-\sum_{j \in J} c_j^{fi}\hat{Y}^*_j+Q^{\omega}(\hat{Y}^*)-z_{N^\prime})^2$, and $z_\alpha$ is the $\alpha$-critical value of the standard normal distribution.} 

\textcolor{mycolor}{\textbf{5.}} Report $\hat{Y}^*$ as the $\varepsilon$-optimal solution with optimality gap \textcolor{mycolor}{$\varepsilon=(U-L)/U$}.
\end{algorithm}

\begin{algorithm}[h] 
\LinesNumberedHidden
\spacingset{1.2}
\footnotesize
  \renewcommand{\arraystretch}{0.4}
    \caption{L-shaped decomposition method}
      \label{alg:BD}
      \SetAlgoNoLine
      \DontPrintSemicolon
\KwIn{Set of optimality cuts \textcolor{mycolor}{$\mathcal{L}$} $\gets \emptyset$, termination criteria $\epsilon \gets 10^{-5}$, lower bound $lb \gets -\infty$, and upper bound $ub \gets \infty$.}
\KwOut{Optimal solution $(z_N,\hat{Y})$ of SAA problem \eqref{SAA_p}.} 

\textbf{1.} Solve master problem \eqref{TS-MP}, and set $\hat{Y}$ and \textcolor{mycolor}{$ub$} as the optimal solution and objective value of \eqref{TS-MP}.
    \begin{align} 
        \textcolor{mycolor}{ub =} & \textcolor{mycolor}{\max \limits_{Y,\theta} \left \lbrace -\sum_{j \in J} c_j^{fi}Y_j + \theta \right \rbrace } \label{TS-MP} \\
        \text{s.t.} & \quad \eqref{C1.1}-\eqref{C1.2} \nonumber \\
        & \quad \textcolor{mycolor}{\theta \leq  \frac{1}{N} \sum_{\omega = 1}^N (h^\omega-T^\omega Y)\pi_l^{\omega}  \quad \forall l \in \mathcal{L}} \nonumber
    \end{align}
   
\textbf{2.} For each scenario $\omega$, calculate $p^{o,\omega}$ and use $\hat{Y}$ to solve sub-problem \eqref{TS-SUBP} and get its duals $\pi^{\omega}$.
    \begin{align}
         \textcolor{mycolor}{Q^{\omega}(\hat{Y}) \coloneqq \max} & \textcolor{mycolor}{\sum_{j \in J} \sum_{k \in K} (p_{k}^{s}-c_j^{pr}-c_{jk}^{t2})V_{jk}^{\omega} + p^{o,\omega} \sum_{j \in J} \bigg ( \sum_{k \neq j \in K} V_{jk}^{\omega} + \xi_{j}^{F^{\prime}\mspace{-7mu},\omega}V_{jj}^{\omega} \bigg ) - \sum_{i \in I}\sum_{j \in J}(c_i^{rm}+c_{ij}^{t1})U_{ij}^\omega} \label{TS-SUBP}\\
         \text{s.t.} & \quad \eqref{C2.1}- \textcolor{mycolor}{\eqref{C2.9}} \nonumber
    \end{align}
    
\textbf{3.} Compute $z_N(\hat{Y})= \textcolor{mycolor}{-\sum_{j \in J} c_j^{fi}\hat{Y}_j + \frac{1}{N} \sum_{\omega=1}^{N} Q^\omega(\hat{Y}).}$ 

\textbf{4.} \textcolor{mycolor}{\lIf {$z_N > lb$}{$lb \gets z_N$.}}

\textbf{5.} \lIf {$[(ub - lb) / ub] > \epsilon$}{
    build an optimality cut \textcolor{mycolor}{$\theta \leq \frac{1}{N} \sum_{\omega = 1}^N (h^\omega-T^\omega Y)\pi^{\omega}$, where $h^\omega$ and $T^\omega$ represent the requirements vector and technology matrix of model \eqref{TS-SUBP}, respectively. Add this cut to the set of constraints $\mathcal{L}$, and go to 1.}}
    \lElse {stop, return $lb$ as the optimal objective function value ($z_N$) and $\hat{Y}$ as the optimal solution.}
\end{algorithm}

The SAA method builds optimization problems through an exterior sampling approach. \textcolor{mycolor}{It estimates $\sum_{\omega \in \Omega} \mathcal{P}^\omega Q^{\omega}(Y)$ by a sample average function $\frac{1}{N} \sum_{\omega=1}^{N} Q^\omega(Y)$}, in which the uncertainty set $\Omega$ is replaced by a set of $N$ i.i.d. sampled scenarios \citep{Verweij2003}. Solving the SAA problem \eqref{SAA_p} for $M$ independent replications, each of size $N$, we get $M$ objectives values $(z_N^1,...,z_N^M)$ and candidate solutions $(\hat{Y}^1,...,\hat{Y}^M)$. \textcolor{mycolor}{An approximate $100(1-2\alpha)\%$ confidence of the statistical estimates of the lower $(L)$ and upper $(U)$ bounds} of the optimal objective value of \eqref{2SSP} can be calculated as presented in Algorithm \ref{alg:SAA-BD}. In particular, \textcolor{mycolor}{$U$} is estimated by solving \eqref{SAA_p} with $N^\prime \gg N$ number of scenarios. This procedure generates a valid bound of the true gap with confidence of at least 1-2$\alpha$. Discussions about the properties of the estimators, convergence of the SAA method and the selection of $M,N$ and $N^\prime$ to get $\varepsilon$-optimal solutions are provided in \citet{Kleywegt2002}.

The L-shaped method (Algorithm \ref{alg:BD}), also called Benders' decomposition, divides the two-stage program into a master problem \eqref{TS-MP} and sub-problem \eqref{TS-SUBP} \citep{SANTOSO200596}. The master problem is a relaxation of \eqref{SAA_p} in which \textcolor{mycolor}{$\frac{1}{N} \sum_{\omega=1}^{N} Q^\omega(Y)$} is approximated by a variable $\theta$ using Benders' cuts built iteratively with the solution of $N$ subproblems \eqref{TS-SUBP}. Since the model has complete recourse, we do not generate feasibility cuts; \textcolor{mycolor}{a single Benders' optimality cut is added in each iteration to the set of optimality cuts $\mathcal{L}$ in the form $\theta \leq \frac{1}{N} \sum_{\omega = 1}^N (h^\omega-T^\omega Y)\pi^{\omega}$, where $\pi^\omega, h^\omega$ and $T^\omega$ represent the dual multipliers, requirements vector, and technology matrix of model \eqref{TS-SUBP}, respectively, following the conventions of \citet{Birge_Louveaux_2011}. }

\section{Case study} \label{s:case}
The proposed model and solution method are implemented in a case study based on vincristine, an essential generic injectable oncology drug used to treat leukemia and Hodgkin's and non-Hodgkin's lymphomas. It has been in shortage in several countries, e.g., in the US \citep{Utah2022}, Colombia \citep{SabogalDeLaPava2022}, and Spain \citep{CIMA2022}. 
\textcolor{mycolor}{In our case study, we consider 11 suppliers of raw materials ($I$), 60 countries to locate plants ($J$), and 179 countries with demand ($K$). These were selected based on the availability of public data. We consider the sets of allies $F_k$ for each country $k \in K$ as the countries with Free Trade Agreements with country $k$ \citep{WTO2024}.} Note that the case study is designed to reflect realistic dynamics of practice not to predict specific drug availability by country.

Parameters are estimated from the following secondary sources. Fixed facility $(c_j^{fi})$, production $(c_j^{pr})$, and raw materials $(c_i^{rm})$ costs are derived from \citet{Tucker2020} and global indexes related to manufacturing activities \citep{WorldBank2020,BCG2020}. \textcolor{mycolor}{Baseline drug prices $(p_k^s)$ are estimated based on price information presented by \citet{Martei2020} and \citet{RR-2956-ASPEC}}. For transportation costs $(c^{t1}_{ij};c^{t2}_{jk})$, we use the International Merchandise Trade - Transportation Costs database \citep{UNCTADstat2016}. Reference values for drug exports \textcolor{mycolor}{$(e_k^{F^\prime}; e^{F}_k)$} are constructed using the trade information of \citet{UNComtradeDatabase2019} and the prevalence rate of the three main cancers treated with vincristine \citep{IARC-WHO2020}. The probability of a given country $k\in K$ implementing export bans, \textcolor{mycolor}{$(1-\rho^{\prime}_k)$}, is estimated from World Trade Organization databases (Quantitative Restrictions \cite{WTO2022a} and Trade Monitoring \cite{WTO2022b}). \textcolor{mycolor}{Probability mass functions for capacity strains $(\xi^{ss}_i; \xi^{ps}_j)$ are derived from capacity utilization databases \citep{FRED2022,TradingEconomics2022}}. \textcolor{mycolor}{The Bernoulli distributions that represent capacity disruptions} $(\xi^{sd}_i; \xi^{pd}_j)$ are estimated from \citet{Tucker2020}.
Demand $(\xi^d_k)$ is assumed to be normally distributed, consistent with other studies (e.g., \citet{Zandkarimkhani2020}), and based on incidence rates from the International Agency for Research on Cancer \citep{IARC-WHO2020}. Further details on data are in the supplementary materials. 

We solve the model with $\alpha=1\%$ to build the bounds \textcolor{mycolor}{($L,U$)}, 30 replications ($M$), \textcolor{mycolor}{300} optimization scenarios ($N$), and \textcolor{mycolor}{3000} evaluation scenarios ($N^{\prime}$). \textcolor{mycolor}{\textit{M} is determined based on the number of cores available for parallel computing and reference values from literature, e.g., \cite{Gangammanavar2021}. $N$ is calibrated considering the optimality gap estimator ($\varepsilon$), variance ($\hat{\sigma}_{\bar{z}_{N,M}}^2$), and the consistency of optimal first-stage decisions across replications. The implementation is in Python with Gurobi v10.0.1 and a high-performance computing cluster (30 CPU cores, 46GB RAM). To evaluate result consistency, we perform five base case runs; we observe optimality gaps below 2\% and similar first-stage decisions. Computational time averages 2877s (48 min) with a standard deviation of 171s (2.8 min)}.

To evaluate the potential effects of geopolitical strain on company performance and drug shortages, we conduct the following experiments. We study current conditions in Section \ref{Sec:ResultsBC} and emerging risks of geopolitical strain in Section \ref{Sec:ExpBan}. In Sections \ref{Sec:Price}--\ref{Sec:Alliances}, we consider the effects of policies from external stakeholders seeking to reduce drug shortages. These include price increases, back-shoring, and strategic alliances. \textcolor{mycolor}{In Section \ref{Sec:VSS}, we present the value of including uncertainty in supply chain design. We discuss practical insights in Section \ref{Sec:insights} and present sensitivity analyses in the supplementary materials.}

\subsection{Base case results} \label{Sec:ResultsBC}
First, we present base case conditions using status quo estimates of geopolitical conditions. \textcolor{mycolor}{The average probability of allowing exports of all countries $\big (\frac{1}{|K|}\sum_{k \in K}\rho_{k}^{\prime} \big )$ is 0.977}, and geopolitical strain may occur when the \textcolor{mycolor}{weighted} average global capacity to produce raw materials is lower than $\tilde{r}=80\%$. 

\begin{figure}[h]
   \centering
   \includegraphics[width=10.5cm]{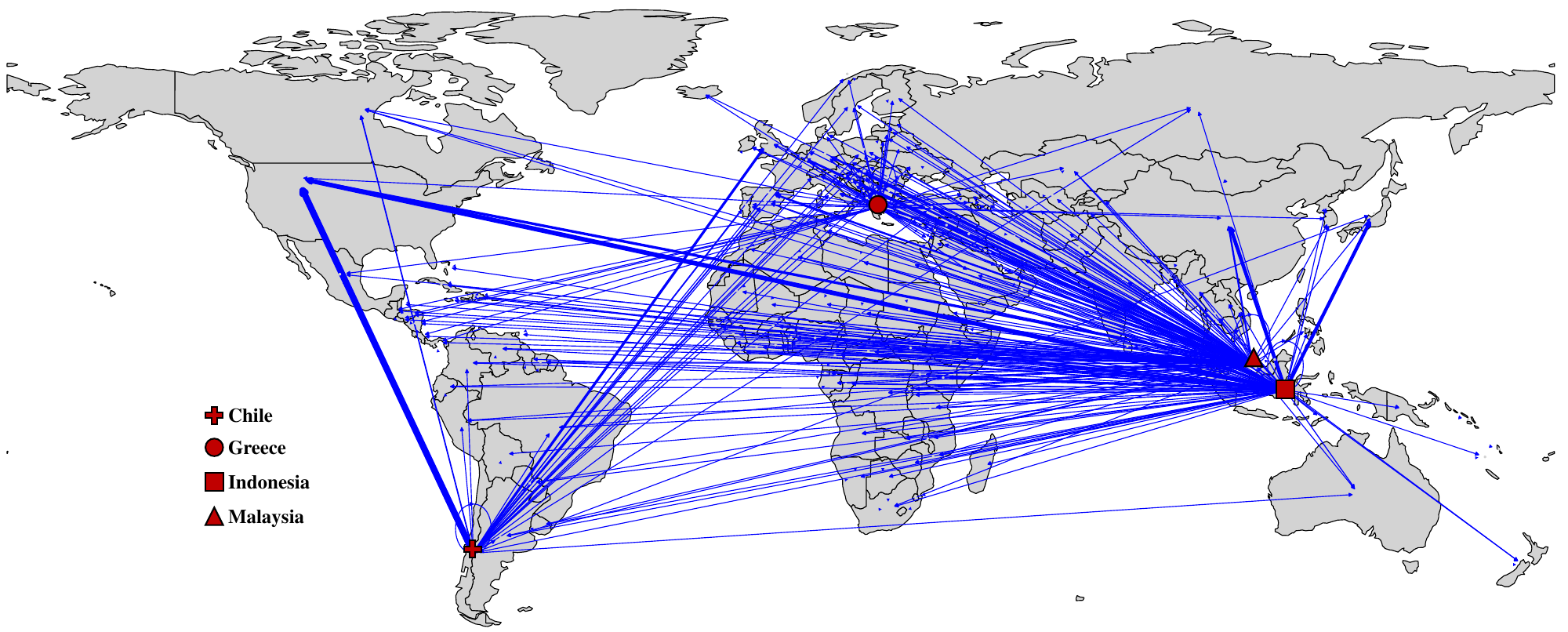}
   \caption{\textcolor{mycolor}{Base case - selected countries and drug distribution.}}
   \label{fig:flowsV-BC.eps}
    \vspace{-10pt}
\end{figure}

Under these conditions, the company locates manufacturing plants in \textcolor{mycolor}{four} countries. Two are countries in Asia (Indonesia and Malaysia), \textcolor{mycolor}{one is in Europe (Greece)} and one is in South America (Chile). Indonesia is a lower-middle income country (LMIC); Malaysia is an upper-middle income country (UMIC); and \textcolor{mycolor}{Greece and Chile are high income countries (HICs)}. No low income country (LIC) is chosen. The selected countries have the lowest fixed and production costs; \textcolor{mycolor}{their probabilities of allowing exports are high, for Indonesia, Malaysia and Chile ($\rho^\prime=0.99$) and Greece ($\rho^\prime=0.9704$)}. The expected available production capacity when there is no disruption $(\mathbb{E}[\xi_{j}^{ps}])$ is 94.9\% (Asian countries), \textcolor{mycolor}{89.0\% (Chile) and 97.8\% (Greece)}. They have low exogenous exports of vincristine. \textcolor{mycolor}{Chile is the third country with the largest number of allies across all countries in set $J$, and Greece does not have allies in our case study}. A map with the selected countries and the distribution of drugs worldwide is in Figure \ref{fig:flowsV-BC.eps}. Note that the greater the intensity and thickness of the arrows, the greater the flow.

Indonesia produces most of the drugs sent to Oceania \textcolor{mycolor}{(97\%)}, Africa \textcolor{mycolor}{(66\%)}, and Latin America \textcolor{mycolor}{(73\%)} and half to Asia \textcolor{mycolor}{(53\%)}. \textcolor{mycolor}{The European volume is produced by Malaysia (35\%), Chile (43\%) and Greece (20\%)}. The North American region (Canada and the US) is attended in \textcolor{mycolor}{51\% by Chile, an ally of Canada and the US, and the 47\% comes from Asian plants}. The top five raw materials suppliers are \textcolor{mycolor}{Thailand, India, Chile, the Czech Republic, and Germany}. The plant located in Chile purchases raw materials from its own country \textcolor{mycolor}{(78\%) and from the US (19\%)}, its ally. The largest opportunity cost \textcolor{mycolor}{($p^o$) is \$11.4} when five countries apply export bans, including the top two drug exporters and all without plant locations. 

\begin{figure}[h]
    \centering
    \vspace{-10pt}
    \includegraphics[width=10cm]{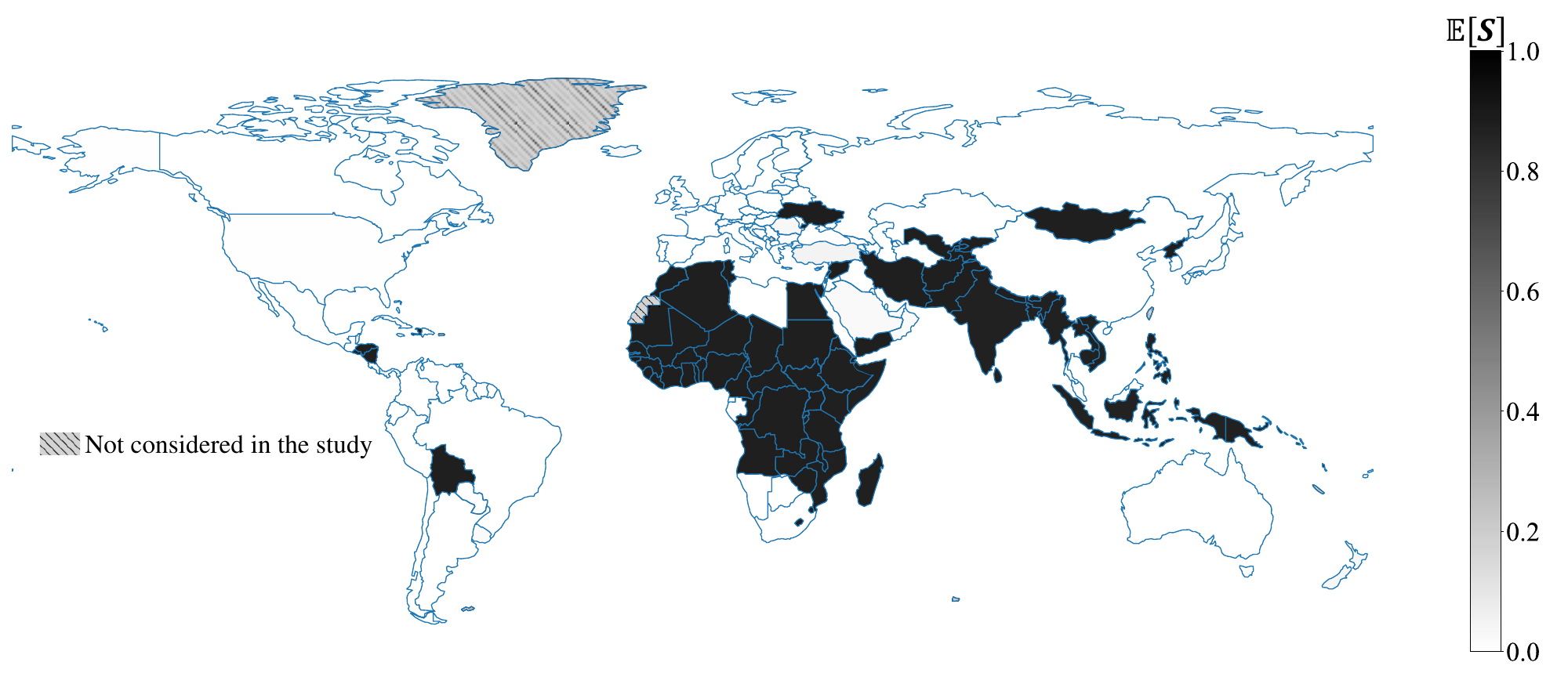}
    \caption{\textcolor{mycolor}{Expected drug shortage $(\mathbb{E}[S])$ per country.}}
    \label{fig:shortage-BC.eps}
    \vspace{-10pt}
\end{figure}

The expected global shortage is \textcolor{mycolor}{17.2\%}, which represents a large portion of unmet worldwide demand. There is an unequal distribution of the drug by region and country income level (Figure \ref{fig:shortage-BC.eps}). The expected shortages for HICs, UMICs, LMICs, and LICs are \textcolor{mycolor}{0.3\%, 0.8\%, 87.2\%, and 87.6\%}, respectively. \textcolor{mycolor}{All LICs and LMICs have an average shortage between 85\% and 90\%}. Notably, one of the countries with a plant, Indonesia (a LMIC), has expected shortages of \textcolor{mycolor}{86.7\%}. This implies that production capacity will not guarantee a domestic drug supply. 

The disparities in drug access persist even in scenarios that have no disruptions. In one scenario without strain, no demand in LMICs and LICs is met despite excess capacity available. \textcolor{mycolor}{In a scenario with a single disruption, demand remains completely met in HICs and UMICs, due to the nominally excess capacity in the fourth plant, and no demand is met in LMICs and LICs. When there are two disruptions, drug shortages are 36\% in HICs and 17\% in UMICs. Sensitivity analyses suggest the results are robust to transportation cost changes.}

\textcolor{mycolor}{We observe that the allies of countries with plant locations have similar expected shortage levels to non-allies. Among LMICs and LICs, alliances are not sufficient to improve access; even allies have expected shortages of around 85\% to 90\%.} 

\subsection{Increasing geopolitical strain} \label{Sec:ExpBan}
\begin{figure}
  \vspace{-10pt}
  \centering
  \includegraphics[width=0.45\textwidth]{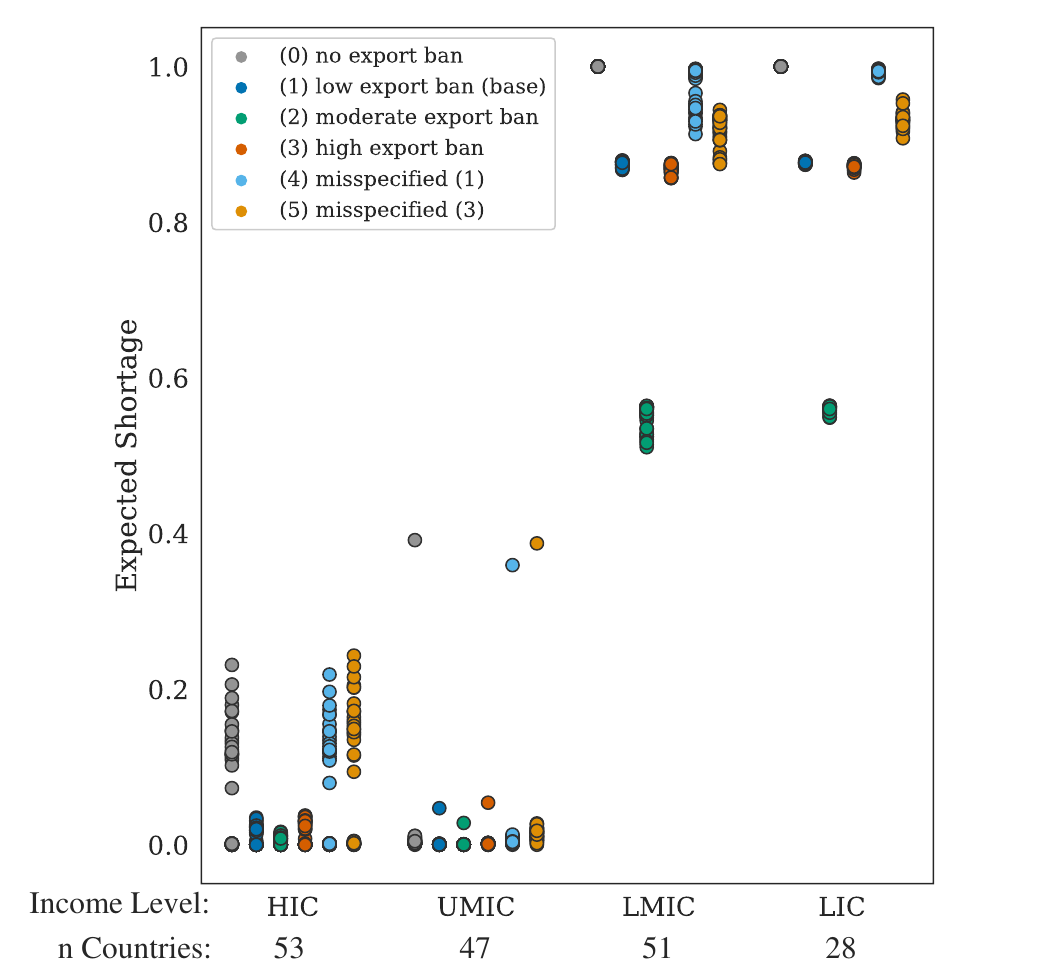}
  \caption{\textcolor{mycolor}{Expected shortage per country under increasing and misspecified geopolitical strain.}}
  \label{fig:known-miss-expban.eps}
  \vspace{-10pt}
\end{figure}

In the wake of the COVID-19 pandemic and increasing international tensions, new geopolitical dynamics could further strain global supply chains. To evaluate potential effects of these emerging risks, we compare six contexts of geopolitical strain. Cases (0) to (3) represent known \textcolor{mycolor}{risks of} none, low (base case), moderate, and high. Cases (4) and (5) represent misspecified plans, where the company selects locations without considering export bans, and the performance is evaluated under low- and high-risk cases. The cases are defined by two parameters: \textcolor{mycolor}{the probability that country $k\in K$ allows exports ($\rho_k^{\prime}$) and the tolerance threshold for raw materials capacity ($\tilde{r}$). Let $\rho^{\prime, base}$ represent the base case values. For country $k\in K$, $(\rho_k^{\prime}, \tilde{r})$ in each case is: case 0 ($1, 80\%$); case 1 ($\rho_k^{\prime, base}, 80\%$); case 2 ($\rho_k^{\prime, base}, 90\%)$; and case 3 $(0.8\rho_k^{\prime, base}, 80\%)$. In cases 4 and 5, the optimization occurs with $\rho_k^{\prime} = 1$}, and the evaluation occurs under low and high-risk conditions, respectively. The risk \textcolor{mycolor}{level} classifications are designated by the expected number of bans each case produces. Figure \ref{fig:known-miss-expban.eps} presents expected drug shortages results in every country under all cases, \textcolor{mycolor}{grouped} them by income level.

\textcolor{mycolor}{First, we consider the known risk of export bans on plant locations. Without export ban risk (case (0)), the selected countries are Indonesia, Malaysia and Greece, a total of three plants. Under low, moderate and high-risk conditions (cases (1) to (3)), four plants are selected. Indonesia and Malaysia are in all cases, and the third and fourth plants are Chile and Greece (case (1)), Greece and Russia (case (2)), and Chile and Colombia (case (3)).} In these configurations, the fourth plant operates as a backup which is used to meet demand under strain; it uses \textcolor{mycolor}{24\%, 34\% and 20\%} of its maximum capacity and produces only \textcolor{mycolor}{8\%, 11\% and 7\%} of worldwide volume, respectively. The selected countries do not have high drug exports, improving market-wide diversification.

\textcolor{mycolor}{Total expected profits increase with increasing export ban risks (7.6\%, 27.5\%, and 32.7\%, vs. case (0) for (1) to (3), respectively). This is due in part to sales volume increases (2\%, 5.1\%, and 1.5\%) and largely due to export ban-induced price increases (70.5\%, 86.3\% and 101.4\% of the revenue increase are attributable to export ban-induced price increases). Costs also increase with increasing export ban risk; fixed costs increase by 34-35\%, and variable costs increase by 5.7\%, 13.1\%, and 2.6\% for each case vs. no risk, respectively.}

\textcolor{mycolor}{We observe that export bans may reduce drug shortages (Table \ref{Table:increasingEB}). Expected global shortages in cases (0) to (3) (none to high risks) are 22.8\%, 17.2\%, 10.7\% and 17.1\%, respectively. The reduction in shortages occurs in part because an additional plant is selected in (1) to (3) that has nominally idle capacity that can be used during disruptions. During export bans, countries are also holding their exports to satisfy domestic demand that may otherwise be unmet; note that the company sales volume increases (6832, 17221, and 5045 units vs. case (0)) are less than the reduction of the global shortage (24081, 51770, and 24193 units vs. case (0)). Expected shortages decrease for countries of all income levels under cases (1) to (3) (Table \ref{Table:increasingEB}), and there are substantial improvements for LMICs and LICs}. This may be because the ban-induced price increases help make satisfying LMIC and LIC demand more attractive. Yet they remain very high and are always worse than in HICs and UMICs.

\begin{table}[h!]
\centering
\renewcommand{\arraystretch}{0.85}
\scriptsize
\caption{\textcolor{mycolor}{Expected shortage by income level under increasing geopolitical strains.}} 
    {\color{mycolor}\begin{tabular}{lcccc}
    \hline
    \textbf{Case}& \textbf{HIC} & \textbf{UMIC} & \textbf{LMIC} & \textbf{LIC}\\ 
    \hline
     (0) no export ban              & 3.6\%  & 6.6\%   & 100\%     & 100\%     \\  \hdashline 
     (1) low export ban (base)      & 0.3\%  & 0.8\%   & 87.2\%    & 87.6\%    \\
     (2) moderate export ban        & 0.2\%  & 0.5\%   & 54\%      & 55.7\%    \\
     (3) high export ban            & 0.5\%  & 0.9\%   & 86.5\%    & 87.1\%    \\
    \hline
    \end{tabular}\label{Table:increasingEB}}
    \vspace{-10pt}
\end{table}

\textcolor{mycolor}{Under the misspecified cases (4) and (5), we observe that when the company does not plan for export ban risks, the company's performance and drug shortages may be negatively impacted. The total expected profits decrease by 1.1\% in the misspecified low-risk case ((4) vs. (1)) and 3.2\% in the misspecified high-risk case ((5) vs. (3)). In both cases, one fewer plant is selected, indicating less resilience. Expected revenue decreases by 4.3\% and 5.3\% with a sales volume reduction of 1.7\% and 2.3\% for low and high-risk, respectively. Expected drug shortages increase substantially when companies do not plan for the risk of export bans; expected shortages globally increase by 29.1\% for low-risk (relative to baseline of 17.2\%) and by 26.4\% for high-risk (relative to baseline of 17.1\%) cases. In the misspecified high-risk case, shortages in HICs, UMICs, LMICs, and LICs increase by 4.1, 6.2, 3.9, and 5.8 absolute percentage points (ppt), respectively vs. properly specified risk. These suggest that under emerging geopolitical threats, ignoring export ban risks during design decisions may decrease the company's expected profits and substantially increase drug shortages.}

\subsection{Pricing policies and impacts on LMICs and LICs} \label{Sec:Price}
\textcolor{mycolor}{Results in Sections \ref{Sec:ResultsBC} and \ref{Sec:ExpBan}} suggest an inequitable distribution of the drug where LMICs and LICs have high shortages. This is consistent with current trends observed in practice \citep{Suda2022} due to the price differentials \citep{Martei2020}. \textcolor{mycolor}{Theorem 2 establishes that economic feasibility is a necessary condition to sell drugs to country $k\in K$, and the case study results indicate} that the low prices in LMICs and LICs are not sufficient to cover the cost of producing and distributing the drugs, even when capacity is available.

To evaluate the potential of pricing policies to promote access in LMICs and LICs, we run three experiments. Each changes the parameter $p_{k}^{s}$, i.e., the baseline \textcolor{mycolor}{drug price in country $k\in K$, and the geopolitical risk is set to the baseline of low export ban risk (case (1)).} In the first, the policy is not to use differential pricing; all countries have the same price (set equal to the US base). In the second and third experiments, we increase LMIC and LIC baseline prices by 50\% and 100\%, respectively. The prices for UMICs and HICs are set to their base values. The latter two represent policies in which governments or non-profit institutions subsidize part of the price, similar to the financial support of GAVI for LIC-LMIC vaccine procurement. Other studies, such as \citet{Kazaz2023} and \citet{DONG-LI2022}, support this idea of social investors' subsidies to supply chains to increase manufacturing capacity or compensate for manufacturing costs when there are product access issues. 

\textcolor{mycolor}{In each experiment, the selected locations are the same (Indonesia, Malaysia, Chile, and Greece). Compared with base case (1), the expected profits reduce by 24.9\% for the no differential prices policy and increase by 3.8\% and 8.4\% for the 50\% and 100\% price supplement policies. The sales volume increases by 22.2\%, 22.2\%, and 22.3\%, respectively. This shows that the reduction in profits for the no differential prices policy is due to the US base price being below prices in other HICs and UMICs, i.e., some HICs and UMICs prices are reduced when set equal to the US price.} 

\textcolor{mycolor}{Global drug access substantially improves in each pricing policy experiment. The expected global drug shortage is 2.4\% for the experiment with equal prices and 2.5\% and 2.6\% with LIC-LMIC price supplements (vs. 17.2\% at base case). The results by country income level are in Table \ref{Table:pricing} and compared with the base case (1) and misspecified low risk case (case (4)). Price supplements are highly effective at reducing shortages for LMICs and LICs; expected shortages for LMICs and LICs drop from baseline values of 87-88\% to under 12\%. There are small increases in shortages for UMICs; this may be because the company switches some of the UMIC supply to LMICs and LICs. The comparator, no differential prices, also shows substantial declines in shortages for LMICs and LICs because US prices are generally higher than LMIC and LIC prices, even with 50\% or 100\% supplements; there are minor increases in shortages for HICs and UMICs. These experiments suggest that policies that increase prices in LMICs and LICs may improve access.}

\setlength{\tabcolsep}{5pt}
\begin{table}[h!]
\centering
\renewcommand{\arraystretch}{0.85}
\caption{Expected shortage by income level under different pricing policies.} 
    \begin{tabular}{lcccc}
    \scriptsize
    \vspace{-10pt}
        {\color{mycolor}\begin{tabular}{lcccc}
        \hline
        \textbf{Experiment}& \textbf{HIC} & \textbf{UMIC} & \textbf{LMIC} & \textbf{LIC}\\ 
        \hline
         Misspecified base (case (4))   & 3.8\%  & 6.1\%   & 96.4\%   & 99.2\%     \\ 
         Base (case (1))                & 0.3\%  & 0.8\%   & 87.2\%   & 87.6\%     \\  \hdashline 
         No differential prices        & 3.3\%  & 1.3\%   & 0.5\%    & 0.7\%      \\
         50\% LIC-LMIC                 & 0.3\%  & 0.9\%   & 10.7\%   & 11.8\%     \\
         100\% LIC-LMIC                & 0.4\%  & 2.8\%   & 8.9\%    & 10.6\%      \\
        \hline
        \end{tabular}\label{Table:pricing}}
    \end{tabular}
\end{table}

\subsection{Back-shoring} \label{Sec:Back}
Next, we study the potential effects of back-shoring, a policy promoted by governments, including the US \citep{TheWhiteHouse2021}. To enforce back-shoring, we require the company to open a plant in the US by adding the constraint $Y_{USA} = 1$ to program \eqref{TS-MP}. We study this under strain cases \textcolor{mycolor}{low-risk} (1) and \textcolor{mycolor}{high-risk} (3). As a secondary analysis, we consider the case where the US plant is of higher quality than base case, i.e., decrease disruption rate ($\xi_{US}^{pd}$) by 50\% (moderate) and 75\% (high) and set operational quality ($\xi_{US}^{ps}$) to the worldwide highest.

When back-shoring to the US, keeping the base case plant quality requirements, the selected countries are Indonesia, Malaysia, Chile, and the US, for both export ban risks. The US manufacturing plant possesses around \textcolor{mycolor}{8\% (low-risk) and 7\% (high-risk)} of the total production, and it uses \textcolor{mycolor}{24\% (low-risk) and 21\% (high-risk) of its maximum manufacturing capacity while the others operate above 85\% (low-risk) and 83\% (high-risk)}. This is because the US plant has higher production costs, and it is only used as a backup to mitigate possible disruptions in the other selected manufacturing plants or the risk of export bans. 

\textcolor{mycolor}{With back-shoring, four manufacturing plants are still selected. Under low-risk, the US replaces Greece (case (1)), and under high-risk the US replaces Colombia (case (3)). The company expected profits reduce by 0.7\% (low-risk) and 0.5\% (high-risk), with increases in fixed costs by  4.2\% (low-risk) and 3.9\% (high-risk). The sales volume reduces by 0.1\%, equivalent to 264 units (low-risk), and 0.3\%, equivalent to 1127 units (high-risk), with increases in variable costs by 1.3\% (low-risk) and 0.7\% (high-risk).} Although there is a plant in the US, Chile supplies the majority of US drugs in expectation, \textcolor{mycolor}{64\%} (low-risk) and \textcolor{mycolor}{63\%} (high-risk). The US supplies \textcolor{mycolor}{5\%} (low-risk) and \textcolor{mycolor}{4\%} (high-risk) of its drugs. Moreover, most of the destinations of the drugs sent by the US are HICs, \textcolor{mycolor}{87.9\%} (low-risk) and \textcolor{mycolor}{86.5\%} (high-risk). These are reflective of the high production cost in the US. 

\textcolor{mycolor}{Drug shortages are fairly similar with back-shoring and properly specified export ban risks. Under low-risk geopolitical strain, with required back-shoring, US shortages are 0.6\% vs. base case of 0.5\%. In a high-risk setting, with required back-shoring or not, US shortages are 0.8\%. US allies do not experience substantial changes when requiring back-shoring to the US.}

In the secondary set of experiments with higher quality US plants, we observe that improving the quality of the production process does not affect the company's economic performance, location decisions, nor global shortages in expectation. \textcolor{mycolor}{There are minor increases in sales volume (between 0.1\% to 0.3\%) that contribute to minor reductions of expected shortages in HICs and UMICs (until 0.6 ppt). For LMICs and LICs, changes (increases and reductions) are below 1 ppt. With the higher quality US plants, there are small improvements to US shortages.}
In the low-risk case, the US shortage was \textcolor{mycolor}{0.4\% (moderate quality increase) and 0.3\% (high quality) vs. 0.6\% back-shoring only}, and in the high-risk case, the US shortage was \textcolor{mycolor}{0.6\% (moderate) and 0.4\% (high) vs. 0.8\% back-shoring only}. 
Shortages in US-allied countries in expectation did not substantially change, even when allied countries have a greater chance of not being affected by export bans when locating a plant in the US. \textcolor{mycolor}{For US-allied HICs and UMICs, under the low-risk case, the average reduction is 0.5 and 0.6 ppt, with the maximum reduction in a country being 1.8 and 2.2 ppt. Under high-risk levels, the average reduction is 0.7 and 0.8 ppt with maximums of 2.7 and 3.4 ppt. For US-allied LMICs, changes in drug shortages are below 1 ppt.} 

\textcolor{mycolor}{These results suggest that back-shoring does not have large impacts on company profitability or shortages when a company is also planning for export ban risks. If a company misspecifies export ban risks, the expected US shortages are 7-9\%, and back-shoring may be worthwhile. }

\subsection{Bilateral alliances} \label{Sec:Alliances}
Alliances have been proposed to reduce drug shortages \citep{TheWhiteHouse2021}. To evaluate the potential effect of bilateral alliances, we perform two experiments in which we remove the \textcolor{mycolor}{effects of alliances by setting the export parameters to the non-allied values}, i.e., \textcolor{mycolor}{$\xi_k^{F} := \xi_k^{F^{\prime}} \forall k \in K$}.

\textcolor{mycolor}{The results suggest that having bilateral alliances to mitigate geopolitical strain does not strongly impact drug shortages; they may affect supply chain decisions with a minor impact on the company's economic performance.}

\textcolor{mycolor}{In the low-risk case without alliances, one of the selected locations (Colombia) is different than in the alliances model (Greece, case (1)), though Greece does not have allies in the case study. Expected profits increase by 0.3\%, fixed costs increase by 0.3\%, and sales volume and variable costs are at the same levels as the alliances model. This suggests that the slight increase in profits is due to higher unit export ban price increases $(p^o)$ since there are more retained exports when there are no alliances. Alliances produce small changes in expected shortages; global shortages are reduced by 0.6\% (0.1 ppt), and shortages in HICs, UMICs, LMICs, and LICs reduce by 0.1, 0.2, 0.1, and 0.3 ppt, respectively. In addition, for the allies of Indonesia, Malaysia and Chile (selected countries in both alliance and non-alliance models), having alliances with countries that host plants results in minor reductions of expected shortages (below 1 ppt, average of 0.2 ppt); six countries experience minor shortage increases (up to 0.4 ppt).} 

\textcolor{mycolor}{In the high-risk case, the same supply chain configurations are selected in both models (with and without alliances). Having the same plants but without the benefits of the alliances reduces expected profits by 1\%, with a decrease in mean sales volume of 0.2\%. Under high-risk of export bans, bilateral alliances improve shortages. Global expected shortages reduce by 4.7\% (0.8 ppt) and shortages in HICs, UMICs, LMICs, and LICs reduce by 0.5, 0.5, 2.1, and 2.5 ppt, respectively. For the allied countries of the hosts of manufacturing plants, their shortages decrease by an average of 1.4 ppt, with a maximum decrease in a country of 3.6 ppt.}

\subsection{\textcolor{mycolor}{Value of considering uncertainty}} \label{Sec:VSS}

\textcolor{mycolor}{We analyze the value of using uncertain parameters. First, we determine the value of the stochastic solution (VSS) under low and high-risk settings. Then, we evaluate the value of planning against the risk of export bans.}

\textcolor{mycolor}{To calculate the VSS, we solve the mean value problem (MVP) of model \eqref{2SSP} and evaluate the MVP optimal solution under the stochastic program \citep{Birge_Louveaux_2011}. The MVP is formulated by setting the stochastic parameters equal to their expected values. Details are in the supplementary materials. The VSS corresponds to the difference between the optimal objective function of the stochastic program \eqref{2SSP} and the objective function value of the MVP solution evaluated in the stochastic program. For the second set of experiments, the proposed stochastic program \eqref{2SSP} is solved setting the export availability parameters, $\xi_{k}^{F^{\prime}}$ and $\xi_{k}^F$, equal to their expected values in the optimization phase and evaluating the first-stage solutions under low and high-risk conditions using the full distributions. To calculate the expectations of $\xi_{k}^{F^{\prime}}$ and $\xi_{k}^F$, we condition on the trigger for their occurrence, i.e., the weighted average global capacity for raw materials, $\big(\sum_{i \in I}q_{i}^{sc}\xi_{i}^{sc}/\sum_{i \in I}q_{i}^{sc}\big) < \tilde{r}$; see supplementary materials for details.}

\textcolor{mycolor}{For low export ban risks, the VSS is \$7841 (0.4\%) and for high export ban risks, the VSS is \$6687 (0.3\%). Supply chain planning considering the mean value of the stochastic parameters would increase global expected shortages by 27.5\% (low-risk) and 0.7\% (high-risk), equivalent to 4.7 ppt (low-risk) and 0.1 ppt (high-risk). Considering country income level classifications, under low export ban risks, HICs, UMICs, LMICs and LICs shortages increase by 3.2, 5.1, 8.9, and 11.4 ppt, respectively. Under high-risk levels, shortages do not change for HICs; however, increase for UMICs, LMICs and LICs by 0.1, 0.4, 0.7 ppt, respectively. The substantial change in shortage under low-risk conditions is because the MVP solution is formed by three manufacturing plants (instead of four), affecting supply chain resilience. This implies that when bans are misspecified, the patients may be the ones predominantly affected, rather than the company itself.}

\textcolor{mycolor}{Results related to planning considering the expected risk of export bans show that if the export ban risk is low, planning using the expected value of the risk of the export ban does not substantially affect the company's performance and drug shortages; using ($\mathbb{E}[\xi_k^{F}];\mathbb{E}[\xi_k^{F^{\prime}}]$) or ($\xi_k^{F};\xi_k^{F^{\prime}}$) produces the same supply chain configuration. However, under the high-risk, Greece is selected instead of Colombia (case (3)); the other three plants remain the same. This supply chain structure reduces the expected profits by 0.9\% (\$18559), the fixed cost reduces by 0.3\%, and the sales volume reduces by 0.2\% (843 units). There are also some minor changes in shortage levels. The expected global shortage increases slightly, moving from 17.1\% (case (3)) to 17.3\% (case $\mathbb{E}[\xi_k^{F}];\mathbb{E}[\xi_k^{F^{\prime}}]$). The expected shortage in HICs does not change, while in UMICs, LMICs, and LICs, shortages increase by 0.1, 0.7 and 1 ppt, respectively.}

\subsection{Practical insights} \label{Sec:insights}
Several recent events have increased geopolitical strain, e.g., the export bans implemented during the COVID-19 pandemic and the Ukraine-Russia war. This work studies the effects of such restrictions on supply chain design and drug shortages.

Our results surprisingly suggest that \textcolor{mycolor}{export ban risks may provide sufficient stress on the system to improve resilience.} In the short-term, export bans may reduce shortages in the country imposing the ban and could reduce shortages \textcolor{mycolor}{globally}.
Different export ban risks may produce different supply chain \textcolor{mycolor}{configurations}. Companies may not be planning for geopolitical strain; our findings indicate that the designed supply chain under misspecified export ban risk underperforms when there are risks of export bans. \textcolor{mycolor}{Even planning for the expected values of export bans is worthwhile, though, in high-risk settings, the use of the full distribution would further reduce shortages.} This highlights a need for companies to evaluate risk levels during planning. 

Pharmaceutical companies are for-profit organizations. In search of efficient operations, e.g., \textcolor{mycolor}{improving profits}, they can make decisions that produce inequities between markets. In our results, LICs and LMICs have high drug shortages in all geopolitical risk cases considered, in contrast to HICs and UMICs, whose shortages are relatively low. We show that with pricing policies, access to drugs in LMICs and LICs could be substantially improved. \textcolor{mycolor}{There are also some minor negative side effects on HICs and UMICs shortages}.

We evaluate the policy of back-shoring manufacturing activities to the US. We find that this strategy is not optimal \textcolor{mycolor}{nor does it improve drug shortages over planning for export bans}, but it does not have substantial negative effects. Government incentives could be used to make back-shoring economically attractive for companies. The US plant is used as a backup rather than primary plant because of its high production costs.

Our study includes \textcolor{mycolor}{bilateral alliances between countries to investigate their potential to mitigate geopolitical strain. These alliances show only minor improvements in the company's economic performance and shortages in expectation. It may produce changes in location decisions.} Given the current long and global supply chain structures, it is possible that unconsidered multi-country alliances may be more effective.

\textcolor{mycolor}{We demonstrate the value of considering uncertainty in global pharmaceutical supply chain design. Companies benefit with slightly higher profits, and countries benefit substantially from lower shortage levels. Our results imply that patients are the ones who pay the consequences of companies failing to plan against uncertainty.}

The focus of the model is to give insights into the dynamics of practice but not actual shortages and behaviors. Limitations include the estimation of data through publicly available sources (see case study description and supplementary materials) that may not precisely reflect confidential information related to demand, costs, and capacities. We assume that capacity disruptions are independent between countries. Export bans do not include retaliations. Competition is only implicitly included via exogenous export parameters and market loss with retained exports. The model does not consider economic features of international trade to focus on the geopolitical dynamics. Finally, drug shortages are reported in the context of demand met by a company optimizing its \textcolor{mycolor}{profits}; it may also be met with other mechanisms.

\section{Conclusion} \label{s:conclusion}
Globalization has created interdependencies and increased supply chain vulnerabilities to emerging geopolitical risks. To address unstudied questions of geopolitical strain, this work proposes a two-stage stochastic program to study the global supply chain design problem under export bans in the pharmaceutical industry. We develop \textcolor{mycolor}{structural insights that present conditions under which demand will not be met and how supply is allocated between countries. The practical effects of the economic feasibility conditions are explored in the numerical study; we observe that the low prices in LMICs and LICs may drive disparities and that supplemental pricing policies may improve access. Broadly, we} focus on determining how geopolitical strain, pricing, back-shoring policies, and bilateral alliances between nations may affect the company's economic performance, design decisions, and shortages globally and by income level classification. \textcolor{mycolor}{We observe that export ban risks may provide sufficient stress on the system to enhance resilience and that planning for export risks is important for company profit and drug access.} 

This model provides a foundation for the inclusion of geopolitical strain in pharmaceutical supply chain design. It also calls attention to studying the effects of geopolitical strain on global supply chains in other industries. Future work could incorporate other international economic features of global trade, such as taxes, exchange rates, and transfer prices. In addition, export bans are one way to restrict international trading; future work could be related to integrating other strains, such as export quotas, licensing, and duties.

\section*{Data availability statement} \label{s:datasharing}
The authors confirm that the data sources and descriptions of estimations are available within the article and its supplementary materials. In particular, during pre-acceptance of the article, the complete database is available upon reasonable request, and upon acceptance, it will be made publicly available via GitHub.

\bibliographystyle{apalike2} 
\spacingset{1}
\bibliography{main-ref}

\newpage
\setcounter{section}{0}
\setcounter{figure}{0}
\setcounter{table}{0}
\setcounter{page}{1}

\input{Supp}

\spacingset{1}

\bibliographystyleSupp{apalike2} 
\bibliographySupp{supp-ref}

\end{document}

%% file: Supp.tex
\begin{center}
     \Large \bf {Supplementary Materials for ``Effects of Geopolitical Strain on Global Pharmaceutical Supply Chain Design and Drug Shortages"}
     \end{center}
       
\begin{center}
    Martha L. Sabogal De La Pava and Emily L. Tucker
    \end{center}
    
\spacingset{1.5}

\textcolor{mycolor}{Details on data used for stochastic and deterministic parameters of the model are presented in Sections \ref{SP} and \ref{DP}, respectively. Sensitivity analyses and additional results are in Section \ref{SA}. Finally, in Section \ref{MVP}, details of the mean value problem parameters are presented.}

\textcolor{mycolor}{The case study was developed using publicly available data, including the selection of the model's sets. Supplier countries were chosen using Volza's Catharanthus Roseus Buyers and Suppliers List \citepSupp{s-Volza2022}, as Vincristine is a vinca alkaloid obtained from Catharanthus Roseus \citepSupp{s-Dhyani2022}. Potential plant locations were selected based on data from the Global Manufacturing Cost-Competitiveness Index \citepSupp{s-BCG2020} and the Doing Business score \citepSupp{s-WorldBank2020}. Demand countries were selected considering countries with records in the International Agency for Research on Cancer databases \citepSupp{s-IARC-WHO2020}. As a result, 11 suppliers, 60 potential plant locations, and 179 demand countries were selected. Stochastic and deterministic parameters were estimated for these sets.}

\section{Stochastic parameters} \label{SP}

\subsection{Demand}
Demand in every country $(\xi_k^{d})$ was assumed to be normally distributed, consistent with other works, e.g., \citeSupp{s-Zandkarimkhani2020}. The mean was estimated as follows: first, we got the total number of new people with cancer per country using: \emph{(i)} the cancer incidence rate of Hodgkin lymphoma, Leukemia, and Non-Hodgkin lymphoma, some of the cancers treated with vincristine \citepSupp{s-IARC-WHO2020}; and \emph{(ii)} the population in every country from \citeSupp{s-Worldometer2020}, both databases from 2020. Second, we estimated the amount of drug required per person (ml/person) using the US vincristine demand presented by \citeSupp{s-Tucker2020} and people in the US with the relevant cancers calculated in the previous step. Finally, since cancer treatments are not completely available in lower-income countries \citepSupp{s-Dhyani2022}, the demand of developing countries was adjusted accordingly with the information presented in \citeSupp{s-Suda2022} for the Antineoplastic and Immunomodulators drug category. To calculate the standard deviations, we used the coefficient of variation for developed and developing countries calculated from \citeSupp{s-Suda2022} under the same drug category and the previously estimated means. Here, we assumed that all patients require the same amount of drug for their treatment.

\subsection{Export bans}
We assumed \textcolor{mycolor}{$\xi_k^{F^{\prime}}$ $\sim {\rm Bernoulli}(\rho^{\prime}_k)$ and $\xi_k^{F}$ $\sim {\rm Bernoulli}(\rho_k)$ for all $k \in K$; where $\rho^{\prime}_k$ and $\rho_k$} represent the probability of allowing exports to non-allies and allies, respectively. We estimated \textcolor{mycolor}{(1 - $\rho^{\prime}_k$)} using the Quantitative Restrictions database from 2012-2022 \citepSupp{s-WTO2022a} and the Trade Monitoring database \citepSupp{s-WTO2022b}, considering the export prohibitions of pharmaceutical products, and removing restrictions based on special agreements such as the Chemical Weapons Convention, Convention on International Trade in Endangered Species of Wild Fauna and Flora, among others. Regarding \textcolor{mycolor}{$\rho_k$}, we assumed \textcolor{mycolor}{$\rho_k$} = 0.85 for countries in the databases. For countries with no records on the databases, that is, countries that have not implemented export bans, we used  \textcolor{mycolor}{$\rho^\prime_k$} = 0.99 and \textcolor{mycolor}{$\rho_k$} = 0.9. \textcolor{mycolor}{Note that we sample $\xi_k^{F}$ $\sim {\rm Bernoulli}(\rho_k)$ when $\xi_{k}^{F^\prime\mspace{-7mu},\omega}=0$.} 

\subsection{Capacity strains}
$\xi_i^{ss}$ and $\xi_j^{ps}$ represent strains that affect manufacturing processes for suppliers $i \in I$ and plants $j \in J$. \textcolor{mycolor}{We calculated probability mass functions, $p(n)$ for $n \in$ $\{0.7, 0.75, 0.8, 0.85,...,1\}$}, using histograms based on the capacity utilization databases \citepSupp{s-FRED2022,s-TradingEconomics2022}. For countries with information in the databases, specific distributions were obtained (USA, CAN, ZAF, NGA, IND, CHN, BRA, ARG, COL, AUS, and NZL); for the others, an average per continent was used. We applied separate values for North, Central, and South America.

\subsection{Disruptions}
We assumed $\xi_i^{sd}$ $\sim {\rm Bernoulli}(\rho_i^{sd})$ for all $i \in I$ and $\xi_j^{pd}$ $\sim {\rm Bernoulli}(\rho_j^{pd})$ for all $j \in J$; where $\rho_i^{sd}$ and $\rho_j^{pd}$ represent the probability of having facility $i$ and $j$ available, respectively. These probabilities were estimated using the average time to disruption and recovery presented by \citeSupp{s-Tucker2020}. The same values were used for every facility.

\section{Deterministic parameters}\label{DP}
All costs were adjusted to the year 2021 using the consumer price index for prescription drugs in the US published on the Statista website \citepSupp{s-Statista2022}.

\subsection{Costs for opening and operating manufacturing plants (fixed costs) and production costs}
For 25 countries, we estimated the fixed cost $(c_{j}^{fi})$ and production costs $(c_{j}^{pr})$ of manufacturing plants in each country using the US costs presented by \citeSupp{s-Tucker2020} as a baseline and the Global Manufacturing Cost-Competitiveness Index of 2019 \citepSupp{s-BCG2020}. We calculated the performance of these countries relative to the US index. For the other countries, we used the Ease of Doing Business score of 2019 \citepSupp{s-WorldBank2020}, and the comparison was made using the index of a referee country in every continent from the previous 25.

\subsection{Raw material cost}
We used the same cost for all suppliers, taking the raw material cost of vincristine from \citeSupp{s-Tucker2020}.

\subsection{Transportation costs}
Transportation costs $(c_{ij}^{t1},c_{jk}^{t2})$ were estimated using the International Merchandise Trade - Transportation Costs data published by the United Nations Conference on Trade and Development \citepSupp{s-UNCTADstat2016}. In particular, we took the CIF-FOB margins per kg and 10000 km for all modes of transport and considered drugs as a commodity. Ten databases were available with cost information for the year 2016. We averaged across the ten databases to estimate an average transportation cost per ml per 10000 km for the origins and destinations of eight regions (Africa, Asia, North America, Caribbean, Central America, South America, Oceania, and Europe). Using the calculated transportation costs per unit distance between the respective regions (\$/ml 10k km) and the distance between major cities of every country, we calculated the unit transportation cost between countries in \$/ml. Since we used CIF-FOB margins, if the origin and destination are the same countries, the cost is \$0.

\subsection{\textcolor{mycolor}{Drug prices}}
\textcolor{mycolor}{Baseline drug prices ($p_k^s$) were estimated as follows}. For OECD countries, we used the International Prescription Drug Price Comparisons for unbranded generics without biologics \citepSupp{s-RR-2956-ASPEC} and the US price presented by \citeSupp{s-Tucker2020} as a baseline. For the rest of the countries, we used the median generic prices for vincristine presented by \citeSupp{s-Martei2020} per country income level classification.

\textcolor{mycolor}{The price increase ($p^o$)} was assumed to be a linear function of the exogenous retained exports produced by export bans. The slope of the linear function was 0.00003, estimated considering as $\Delta$prices 50\% and $\Delta$shortages the US mean demand 90000 ml. The FDA has used 50\% as a referent value to analyze shortages of drugs that experience sustained price increases after a shortage occurs \citepSupp{s-FDA2019}. 

\subsection{Production capacity}
Manufacturing plants' and suppliers' capacities $(q_{j}^{pc},q_{i}^{sc})$ were assumed to be 120000 ml and 100000 ml, respectively.

\subsection{Exports}
Reference values for drug exports of vincristine \textcolor{mycolor}{$(e_k^{F^\prime}$,$e_k^{F})$} by volume and per country were constructed as follows: first, using the \citeSupp{s-UNComtradeDatabase2019} and considering HS code 300449 - medicaments containing alkaloids or derivatives, we calculated the percentage of participation of every country in the international trade. Note that vincristine is a vinca alkaloid \citepSupp{s-Dhyani2022}. Second, we defined the total global trade of vincristine using the global prevalence of the cancers mentioned previously (from the International Agency for Research on Cancer website \citepSupp{s-IARC-WHO2020}) and following the same procedure implemented to get the mean demand. We multiply the participation percentage of each country in the international trade by the total exports \textcolor{mycolor}{to get the total volume of exports of country $k$ $(e_k^{F^\prime}+e_k^{F})$. Finally, we set the proportion of exports to allies and non-allies of each country consistently with the proportion of US exports to allies and non-allies (estimated from \citeSupp{s-UNComtradeDatabase2019}); particularly, we assumed that country $k$ provides 0.1\% of their total export to each allied country}.

\section{Sensitivity analyses and additional results}\label{SA}
\subsection{Sensitivity analysis of transportation costs}\label{SA:tc}
Transportation costs are the largest component of the logistics costs of an enterprise \citepSupp{s-Rodrigue2020}. To evaluate the effects of higher transportation costs on the performance and decisions of the company, we ran an experiment with doubled values for $c_{ij}^{t1}$ and $c_{jk}^{t2}$ under the base case conditions. \textcolor{mycolor}{As a result, three location decisions were the same (Indonesia, Malaysia, and Chile), and one was different (Thailand was selected instead of Greece). This supply configuration has the same fixed cost levels, with Thailand's capacity utilization and production levels mirroring those of Greece. The expected profits increased by 0.2\%, and the sales volume remained the same, with a decrease in variable costs by 0.6\%. The total income increased by 0.3\%. The expected global shortage decreased by 0.3 ppt (moving from 17.2\% to 16.9\%). Regarding the expected shortage per income level, for HICs, it did not change; for UMICs, the shortage increased by 0.1 ppt; and for LMICs and LICs, the shortage decreased by 1.8 ppt for both.} These show that our results are robust to changes in transportation costs.

\subsection{Pricing policies additional results}\label{SA:price}

\begin{figure}[H]
   \centering
   \includegraphics[width=7cm]{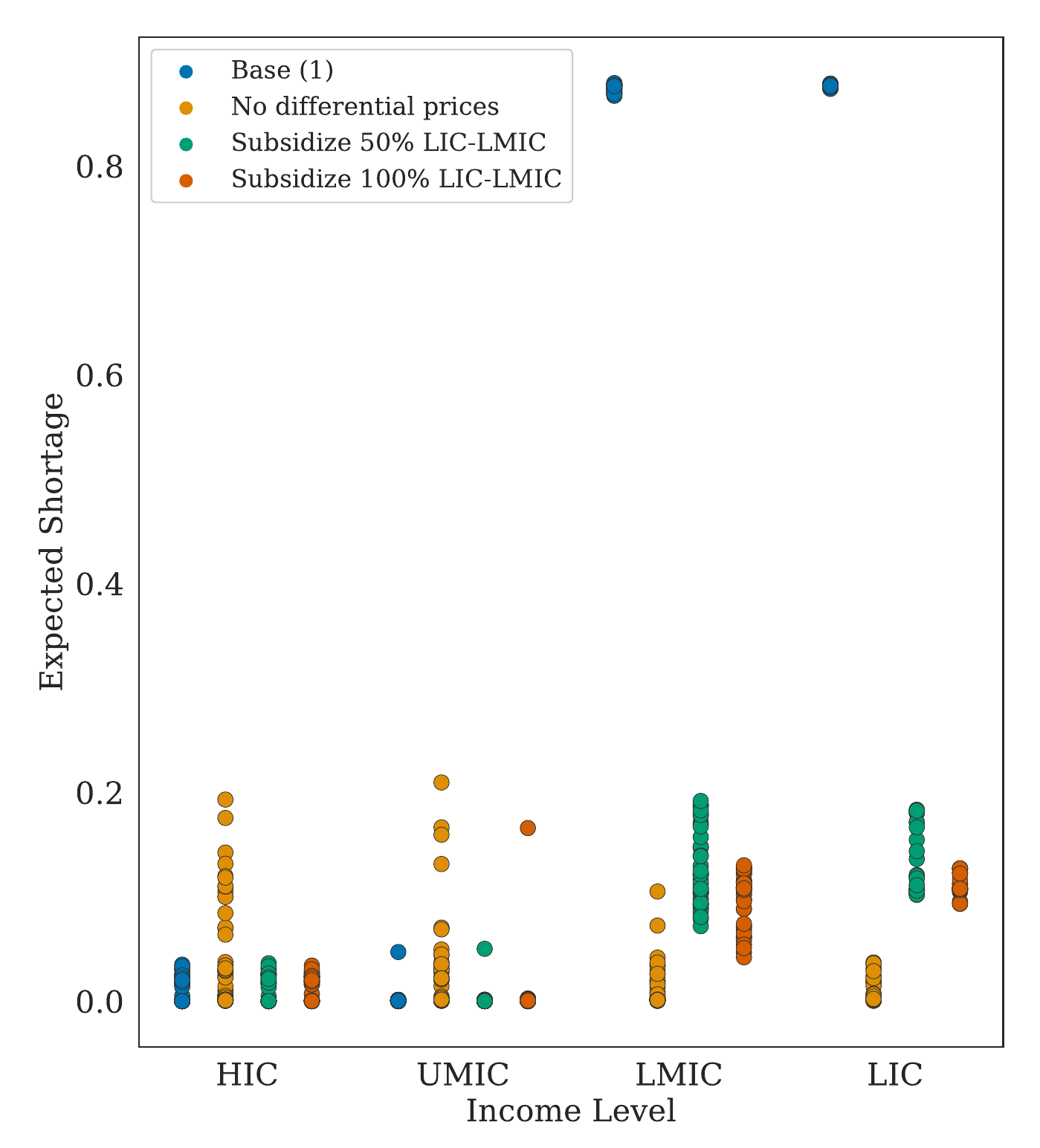}
   \caption{\textcolor{mycolor}{Expected drug shortage under different pricing schemes.}}
   \label{fig:prices.eps}
\end{figure}

\begin{figure}[H]
   \centering
   \includegraphics[width=14cm]{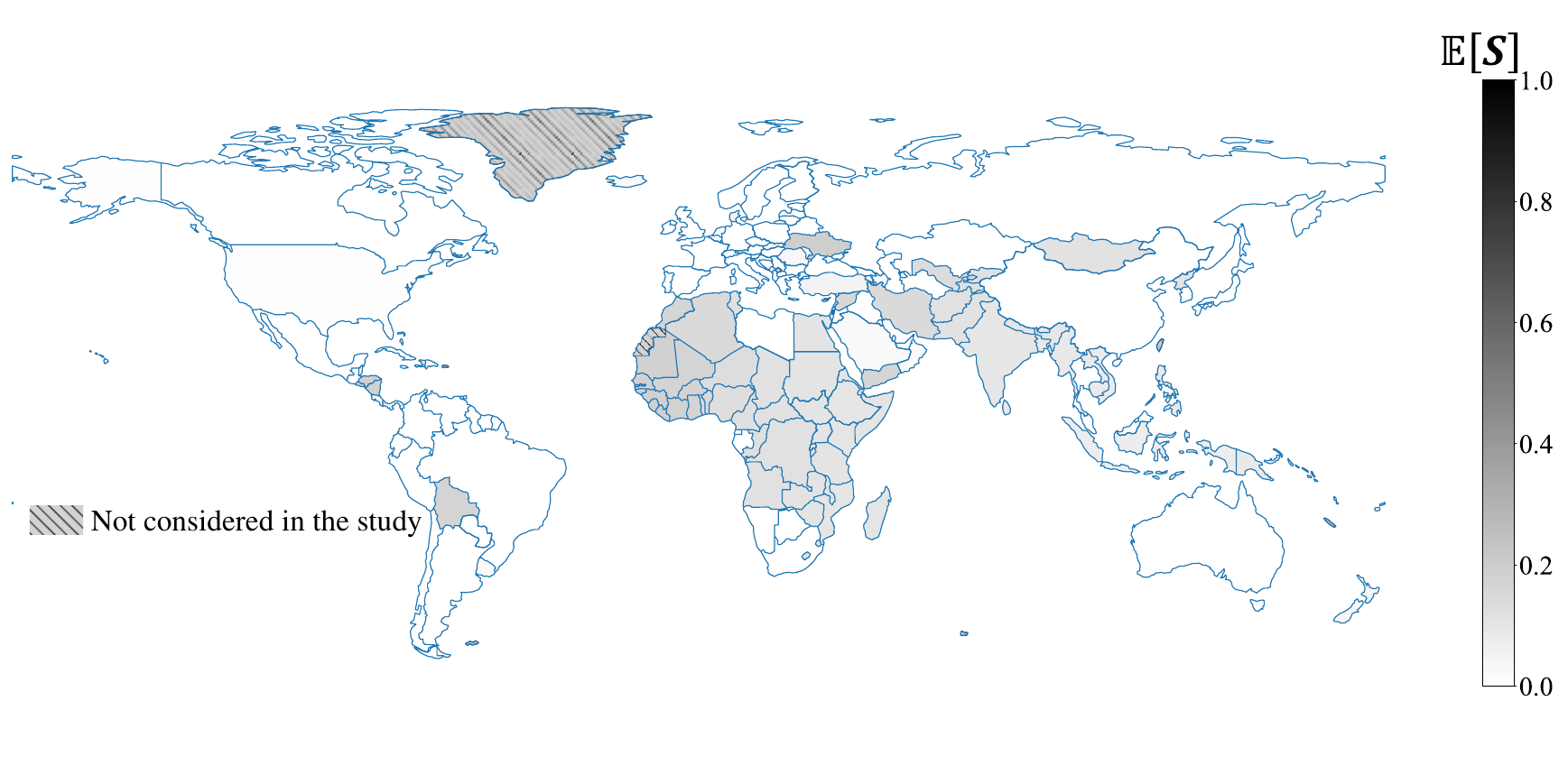}
   \caption{\textcolor{mycolor}{Expected drug shortage $(\mathbb{E}[S])$ under subsidize 50\% LIC-LMIC price policy.}}
   \label{fig:shortage-price50.eps}
\end{figure}

The effect of the studied pricing schemes (main text Section \ref{Sec:Price}) on the expected drug shortage per country grouped by income level is presented in Figure \ref{fig:prices.eps}. Figure \ref{fig:shortage-price50.eps} presents results for the pricing policy subsidizing 50\% LIC-LMIC prices. Particularly, it shows the expected shortage per country.

\vspace{-10pt}

\subsection{Back-shoring additional results}\label{SA:back}

\begin{figure}[H]
    \centering
    \includegraphics[width=14cm]{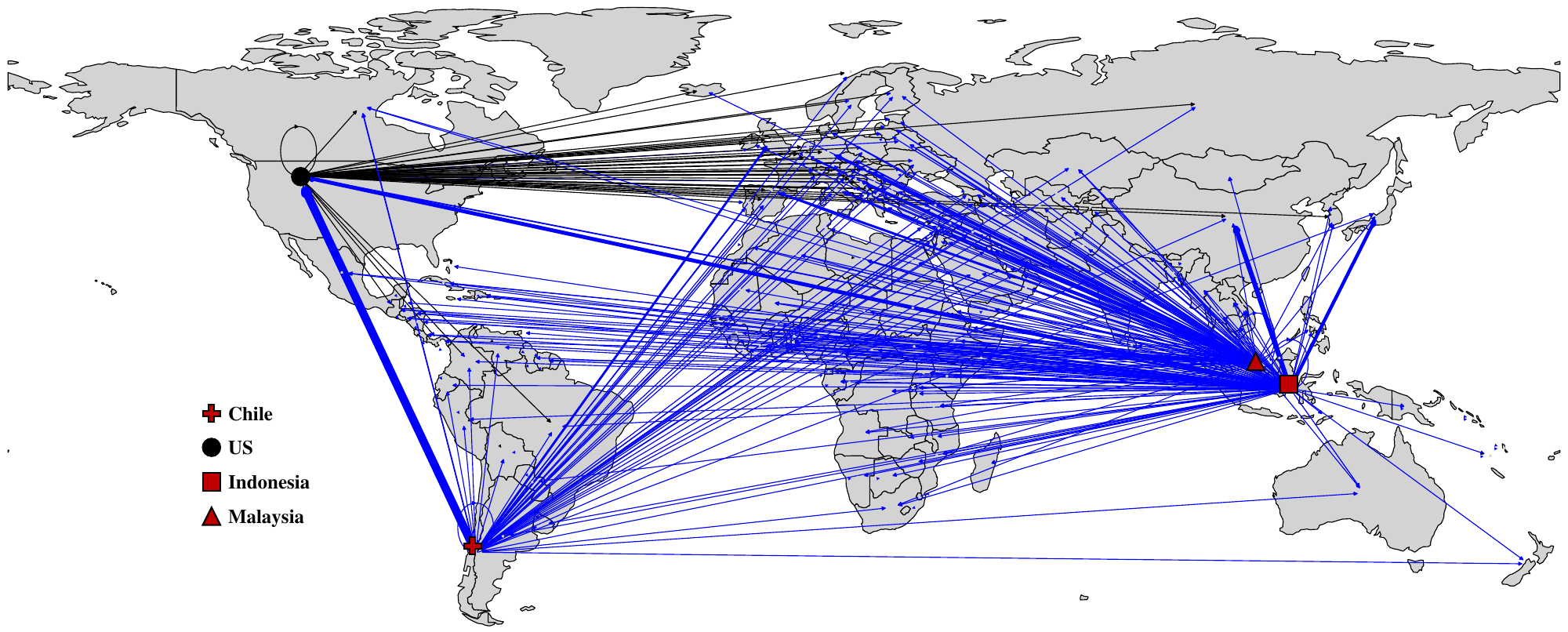}
    \caption{\textcolor{mycolor}{Back-shoring - selected countries and drug distribution under the base case conditions.}}
    \label{fig:backS.eps}
\end{figure}

In the main text Section \ref{Sec:Back}, we present results when the back-shoring policy is required under both geopolitical strain conditions (low- and high-risk) and under three manufacturing plant quality levels (base case, moderate and high quality). Data used in the experiments and additional results are presented in this supplementary section. For the moderate and high-quality experiments, we improved the uncertainty parameters $\xi_{US}^{pd}$ and $\xi_{US}^{ps}$ simultaneously. Specifically, for $\xi_{US}^{pd}$ we decreased the rate of disruptions by 50\% (moderate-quality) and 75\% (high-quality), resulting in probabilities of internal disruptions $\rho_{US}^{pd}=0.9861$ and $\rho_{US}^{pd}=0.9930$, respectively. In addition, for $\xi_{US}^{ps}$ we improved the expected value of the strains that affect the production capacity of the US, i.e., $\mathbb{E}[\xi_{US}^{ps}]$ was 0.8876; in this experiment, it was increased to 0.9868, to be equal to the worldwide highest. Figure \ref{fig:backS.eps} shows the selected locations and the drug distribution worldwide of the case low geopolitical strain and base case quality. 

\section{\textcolor{mycolor}{Parameters of the mean value problem (MVP)}}\label{MVP}
\textcolor{mycolor}{The MVP requires the stochastic parameters of the model to be equal to their expected values. In our model, the stochastic parameters are demand, the availability to export, and the capacity of suppliers of raw materials and manufacturing plants. Descriptions of their distributions are in Section \ref{SP} of this supplementary material. Next, we present how we calculate the expected values of each of these stochastic parameters.}

\textcolor{mycolor}{\textit{Demand ($\xi^d_k$)}. The $\mathbb{E}[\xi^d_k]$ corresponds to the mean of the normal distribution of $\xi^d_k$.}

\textcolor{mycolor}{\textit{Capacity of suppliers $(\xi^{sc}_i)$}. The $\mathbb{E}[\xi^{sc}_i]$ corresponds to the product of the expected value of two independent discrete random variables, i.e., $\mathbb{E}[\xi^{sc}_i] =\mathbb{E}[\xi^{ss}_i]\mathbb{E}[\xi^{sd}_i]$.}

\textcolor{mycolor}{\textit{Capacity of manufacturing plants $(\xi^{pc}_j)$}. The same procedure as the capacity of suppliers.}

\textcolor{mycolor}{\textit{Availability to export $(\xi^{F^\prime}_k;\xi^{F}_k)$}. Both random variables are conditional on the weighted average global availability of raw materials. To simplify notation, let $B$, be the event $\sum_{i \in I}q_{i}^{sc}\xi_{i}^{sc}/\sum_{i \in I}q_{i}^{sc} < \Tilde{r}$, and $B^\prime$ be the complement of $B$. The $\mathbb{E}[\xi^{F^{\prime}}_k] =\mathbb{E}[\xi^{F^{\prime}}_k|B]\mathbb{P}(B)+\mathbb{E}[\xi^{F^{\prime}}_k|B^\prime]\mathbb{P}(B^\prime)=\rho_k^{\prime}\mathbb{P}(B)+1(1-\mathbb{P}(B))$, where $\mathbb{P}(B) = 0.14756$, estimated with an empirical sampling approach, using 50000 samples of $\sum_{i \in I}q_{i}^{sc}\xi_{i}^{sc}/\sum_{i \in I}q_{i}^{sc}$ generated independently. Similar steps are followed to calculate $\mathbb{E}[\xi^{F}_k]$.}

%% file: main_document.bbl
\begin{thebibliography}{}

\bibitem[BCG, 2020]{s-BCG2020}
BCG (2020).
\newblock {A manufacturing strategy built for trade instability}.
\newblock Retrieved from \url{https://www.bcg.com/publications/2020/manufacturing-strategy-built-trade-instability}. Accessed February 16, 2023.

\bibitem[{Board of Governors of the Federal Reserve System US}, 2022]{s-FRED2022}
{Board of Governors of the Federal Reserve System US} (2022).
\newblock {Capacity utilization: Manufacturing: Non-durable goods: Chemical NAICS 325}.
\newblock Retrieved from \url{https://fred.stlouisfed.org/series/CAPUTLG325S}. Accessed December 10, 2022.

\bibitem[Dhyani et~al., 2022]{s-Dhyani2022}
Dhyani, P., Quispe, C., Sharma, E., Bahukhandi, A., Sati, P., Attri, D.~C., Szopa, A., Sharifi-Rad, J., Docea, A.~O., Mardare, I., Calina, D., \& Cho, W.~C. (2022).
\newblock {Anticancer potential of alkaloids: A key emphasis to colchicine, vinblastine, vincristine, vindesine, vinorelbine and vincamine}.
\newblock {\em Cancer Cell International}, 22(1), 1--20.
\newblock \url{https://doi.org/10.1186/s12935-022-02624-9}.

\bibitem[{FDA}, 2019]{s-FDA2019}
{FDA} (2019).
\newblock {Drug shortages: Root causes and potential solutions}.
\newblock Retrieved from \url{https://www.fda.gov/media/131130/download}. Accessed May 1, 2022.

\bibitem[{IARC}, 2020]{s-IARC-WHO2020}
{IARC} (2020).
\newblock {Cancer over time}.
\newblock Retrieved from \url{https://gco.iarc.fr/overtime/en}. Accessed October 1, 2022.

\bibitem[Martei et~al., 2020]{s-Martei2020}
Martei, Y.~M., Iwamoto, K., Barr, R.~D., Wiernkowski, J.~T., \& Robertson, J. (2020).
\newblock Shortages and price variability of essential cytotoxic medicines for treating children with cancers.
\newblock {\em BMJ Global Health}, 5(11), 13.
\newblock \url{https://doi.org/10.1136/bmjgh-2020-003282}.

\bibitem[Mulcahy et~al., 2021]{s-RR-2956-ASPEC}
Mulcahy, A.~W., Whaley, C.~M., Gizaw, M., Schwam, D., Edenfield, N., \& Becerra-Ornelas, A.~U. (2021).
\newblock {\em International Prescription Drug Price Comparisons: Current Empirical Estimates and Comparisons with Previous Studies}.
\newblock Santa Monica, CA: RAND Corporation.

\bibitem[Rodrigue, 2020]{s-Rodrigue2020}
Rodrigue, J.-P. (2020).
\newblock {\em {Trade, Logistics and Freight Distribution}}.
\newblock London: Routledge, 5th edition.
\newblock \url{https://doi.org/10.4324/9780429346323}.

\bibitem[{Statista}, 2022]{s-Statista2022}
{Statista} (2022).
\newblock {Consumer price index for prescription and nonprescription drugs in the U.S. from 1960 to 2021}.
\newblock Retrieved from \url{https://www.statista.com/statistics/187253}. Accessed February 17, 2023.

\bibitem[Suda et~al., 2022]{s-Suda2022}
Suda, K.~J., Kim, K.~C., Hernandez, I., Gellad, W.~F., Rothenberger, S., Campbell, A., Malliart, L., \& Tadrous, M. (2022).
\newblock {The global impact of COVID-19 on drug purchases: A cross-sectional time series analysis}.
\newblock {\em Journal of the American Pharmacists Association}, 62(3), 15.
\newblock \url{https://doi.org/10.1016/j.japh.2021.12.014}.

\bibitem[{Trading Economics}, 2022]{s-TradingEconomics2022}
{Trading Economics} (2022).
\newblock {Capacity utilization}.
\newblock Retrieved from \url{https://tradingeconomics.com/country-list/capacity-utilization}. Accessed December 10, 2022.

\bibitem[Tucker et~al., 2020]{s-Tucker2020}
Tucker, E.~L., Daskin, M.~S., Sweet, B.~V., \& Hopp, W.~J. (2020).
\newblock {Incentivizing resilient supply chain design to prevent drug shortages: Policy analysis using two- and multi-stage stochastic programs}.
\newblock {\em IISE Transactions}, 52(4), 394--412.
\newblock \url{https://doi.org/10.1080/24725854.2019.1646441}.

\bibitem[{UN Comtrade Database}, 2019]{s-UNComtradeDatabase2019}
{UN Comtrade Database} (2019).
\newblock {Global trade data}.
\newblock Retrieved from \url{https://comtradeplus.un.org/}. Accessed December 10, 2022.

\bibitem[{UNCTADstat}, 2016]{s-UNCTADstat2016}
{UNCTADstat} (2016).
\newblock {International merchandise trade - transport costs}.
\newblock Retrieved from \url{https://unctadstat.unctad.org/wds/ReportFolders/reportFolders.aspx?sCS\_ChosenLang=en}. Accessed February 17, 2023.

\bibitem[{Volza}, 2022]{s-Volza2022}
{Volza} (2022).
\newblock {Catharanthus Roseus Export Import Data, Buyers and Suppliers List}.
\newblock Retrieved from \url{https://www.volza.com/p/catharanthus-roseus/}. Accessed October 1, 2022.

\bibitem[{World Bank}, 2020]{s-WorldBank2020}
{World Bank} (2020).
\newblock {\em {Doing Business 2020}}.
\newblock Washington, DC: {World Bank}.
\newblock \url{https://doi.org/10.1596/978-1-4648-1440-2}.

\bibitem[{Worldometer}, 2020]{s-Worldometer2020}
{Worldometer} (2020).
\newblock {Countries in the world by population 2020}.
\newblock Retrieved from \url{https://www.worldometers.info/world-population/population-by-country/}. Accessed October 1, 2022.

\bibitem[{WTO}, 2022a]{s-WTO2022a}
{WTO} (2022a).
\newblock {Quantitative restrictions database}.
\newblock Retrieved from \url{https://qr.wto.org/en\#/home}. Accessed October 1, 2022.

\bibitem[{WTO}, 2022b]{s-WTO2022b}
{WTO} (2022b).
\newblock {Trade monitoring database}.
\newblock Retrieved from \url{https://tmdb.wto.org/en/explore/goods\#nogo\%0A}. Accessed October 1, 2022.

\bibitem[Zandkarimkhani et~al., 2020]{s-Zandkarimkhani2020}
Zandkarimkhani, S., Mina, H., Biuki, M., \& Govindan, K. (2020).
\newblock {A chance constrained fuzzy goal programming approach for perishable pharmaceutical supply chain network design}.
\newblock {\em Annals of Operations Research}, 295(1), 425--452.
\newblock \url{https://doi.org/10.1007/s10479-020-03677-7}.

\end{thebibliography}


\begin{thebibliography}{}

\bibitem[AAM, 2021]{AAM2021}
AAM (2021).
\newblock {A blueprint for enhancing the security of the US pharmaceutical supply chain}.
\newblock Retrieved from \url{https://accessiblemeds.org/sites/default/files/2020-04/AAM-Blueprint-US-Pharma-Supply-Chain.pdf}. Accessed May 2, 2023.

\bibitem[Aguiar \& Kasper, 2021]{EscoAguiar2021}
Aguiar, E. \& Kasper, E. (2021).
\newblock {Trade flows of parallel imported medicines in Europe}.
\newblock Retrieved from \url{https://affordablemedicines.eu/portfolio-item/position-paper-trade-flows-of-parallel-imported-medicines-2020/}. Accessed May 1, 2023.

\bibitem[Akbarpour et~al., 2020]{Akbarpour2020}
Akbarpour, M., {Ali Torabi}, S., \& Ghavamifar, A. (2020).
\newblock {Designing an integrated pharmaceutical relief chain network under demand uncertainty}.
\newblock {\em Transportation Research Part E: Logistics and Transportation Review}, 136, 101867.
\newblock \url{https://doi.org/10.1016/j.tre.2020.101867}.

\bibitem[Aldrighetti et~al., 2021]{Aldrighetti2021}
Aldrighetti, R., Battini, D., Ivanov, D., \& Zennaro, I. (2021).
\newblock {Costs of resilience and disruptions in supply chain network design models: A review and future research directions}.
\newblock {\em International Journal of Production Economics}, 235, 108103.
\newblock http://doi.org/10.1016/j.ijpe.2021.108103.

\bibitem[Ambulkar et~al., 2015]{Ambulkar2015}
Ambulkar, S., Blackhurst, J., \& Grawe, S. (2015).
\newblock {Firm's resilience to supply chain disruptions: Scale development and empirical examination}.
\newblock {\em Journal of Operations Management}, 33-34, 111--122.
\newblock https://doi.org/10.1016/j.jom.2014.11.002.

\bibitem[Azad \& Hassini, 2019]{AZAD2019}
Azad, N. \& Hassini, E. (2019).
\newblock Recovery strategies from major supply disruptions in single and multiple sourcing networks.
\newblock {\em European Journal of Operational Research}, 275(2), 481--501.
\newblock \url{https://doi.org/10.1016/j.ejor.2018.11.044}.

\bibitem[Barbieri et~al., 2020]{Barbieri2020}
Barbieri, P., Boffelli, A., Elia, S., Fratocchi, L., Kalchschmidt, M., \& Samson, D. (2020).
\newblock {What can we learn about reshoring after Covid-19?}
\newblock {\em Operations Management Research}, 13, 131--136.
\newblock \url{https://doi.org/10.1007/s12063-020-00160-1}.

\bibitem[BCG, 2020]{BCG2020}
BCG (2020).
\newblock {A manufacturing strategy built for trade instability}.
\newblock Retrieved from \url{https://www.bcg.com/publications/2020/manufacturing-strategy-built-trade-instability}. Accessed February 16, 2023.

\bibitem[Birge \& Louveaux, 2011]{Birge_Louveaux_2011}
Birge, J.~R. \& Louveaux, F. (2011).
\newblock {\em Introduction to stochastic programming}.
\newblock Springer New York.

\bibitem[Blossey et~al., 2021]{Blossey2021}
Blossey, G., Hahn, G.~J., \& Koberstein, A. (2021).
\newblock {Managing uncertainty in pharmaceutical supply chains: A structured review}.
\newblock In {\em Proceedings of the Annual Hawaii International Conference on System Sciences}  (pp.\ 1435--1444).
\newblock \url{https://doi.org/10.24251/hicss.2021.173}.

\bibitem[Blossey et~al., 2022]{Blossey2022}
Blossey, G., Hahn, G.~J., \& Koberstein, A. (2022).
\newblock {Planning pharmaceutical manufacturing networks in the light of uncertain production approval times}.
\newblock {\em International Journal of Production Economics}, 244(Oct 2021), 108343.
\newblock \url{https://doi.org/10.1016/j.ijpe.2021.108343}.

\bibitem[{Board of Governors of the Federal Reserve System US}, 2022]{FRED2022}
{Board of Governors of the Federal Reserve System US} (2022).
\newblock {Capacity utilization: Manufacturing: Non-durable goods: Chemical NAICS 325}.
\newblock Retrieved from \url{https://fred.stlouisfed.org/series/CAPUTLG325S}. Accessed December 10, 2022.

\bibitem[Bochenek et~al., 2018]{Bochenek2018}
Bochenek, T., Abilova, V., Alkan, A., {...}, \& Godman, B. (2018).
\newblock {Systemic measures and legislative and organizational frameworks aimed at preventing or mitigating drug shortages in 28 European and Western Asian countries}.
\newblock {\em Frontiers in Pharmacology}, 8, 942.
\newblock \url{https://doi.org/10.3389/fphar.2017.00942}.

\bibitem[Casey \& Cimino-Isaacs, 2021]{Casey2021}
Casey, C.~A. \& Cimino-Isaacs, C.~D. (2021).
\newblock {Export restrictions in response to the COVID-19 pandemic}.
\newblock Retrieved from \url{https://crsreports.congress.gov/product/details?prodcode=IF11551}. Accessed May 2, 2023.

\bibitem[Chatzoglou et~al., 2018]{Chatzoglou2018}
Chatzoglou, P., Chatzoudes, D., Petrakopoulou, Z., \& Polychrou, E. (2018).
\newblock {Plant location factors: A field research}.
\newblock {\em Opsearch}, 55(3-4), 749--786.
\newblock \url{https://doi.org/10.1007/s12597-018-0341-1}.

\bibitem[CIMA, 2022]{CIMA2022}
CIMA (2022).
\newblock {Shortage actives and supply problems}.
\newblock Retrieved from \url{https://cima.aemps.es/cima/publico/listadesabastecimiento.html}. Accessed May 2, 2023.

\bibitem[Das, 2018]{Das2018}
Das, T. (2018).
\newblock {\em Managing Trust in Strategic Alliances}.
\newblock Research in Strategic Alliances. Charlotte, NC: Information Age Publishing Inc.

\bibitem[Diaz et~al., 2023]{Diaz2022}
Diaz, R., Kolachana, S., \& Gomes, R.~F. (2023).
\newblock A simulation-based logistics assessment framework in global pharmaceutical supply chain networks.
\newblock {\em Journal of the Operational Research Society}, 74(5), 1242--1260.
\newblock \url{https://doi.org/10.1080/01605682.2022.2077661}.

\bibitem[Duran et~al., 2011]{Duran2011}
Duran, S., Gutierrez, M.~A., \& Keskinocak, P. (2011).
\newblock {Pre-positioning of emergency items for CARE international}.
\newblock {\em Interfaces}, 41(3), 223--237.
\newblock \url{https://doi.org/10.1287/inte.1100.0526}.

\bibitem[Dyer, 2019]{Dyerl6086}
Dyer, O. (2019).
\newblock Us paediatric oncologists are forced to prioritise patients for vincristine treatment as supplies run short.
\newblock {\em BMJ}, 367, l6086.
\newblock \url{https://doi.org/10.1136/bmj.l6086}.

\bibitem[ECLAC, 2021]{ECLAC2021}
ECLAC (2021).
\newblock {Plan for self-sufficiency in health matters in Latin America and the Caribbean: Lines of action and proposals}.
\newblock Retrieved from \url{https://hdl.handle.net/11362/47253}. Accessed May 5, 2023.

\bibitem[Franco \& Alfonso-Lizarazo, 2017]{Franco2017}
Franco, C. \& Alfonso-Lizarazo, E. (2017).
\newblock {A structured review of quantitative models of the pharmaceutical supply chain}.
\newblock {\em Complexity}, 2017, 13.
\newblock \url{https://doi.org/10.1155/2017/5297406}.

\bibitem[Gangammanavar et~al., 2021]{Gangammanavar2021}
Gangammanavar, H., Liu, Y., \& Sen, S. (2021).
\newblock {Stochastic decomposition for two-stage stochastic linear programs with random cost coefficients}.
\newblock {\em INFORMS Journal on Computing}, 33(1), 51--71.
\newblock \url{https://doi.org/10.1287/ijoc.2019.0929}.

\bibitem[Ghavamifar et~al., 2018]{Ghavamifar2018}
Ghavamifar, A., Makui, A., \& Taleizadeh, A.~A. (2018).
\newblock {Designing a resilient competitive supply chain network under disruption risks: A real-world application}.
\newblock {\em Transportation Research Part E: Logistics and Transportation Review}, 115(2018), 87--109.
\newblock https://doi.org/10.1016/j.tre.2018.04.014.

\bibitem[Goh et~al., 2007]{Goh2007}
Goh, M., Lim, J.~Y., \& Meng, F. (2007).
\newblock {A stochastic model for risk management in global supply chain networks}.
\newblock {\em European Journal of Operational Research}, 182(1), 164--173.
\newblock \url{https://doi.org/10.1016/j.ejor.2006.08.028}.

\bibitem[Goodarzian et~al., 2020]{Goodarzian2020}
Goodarzian, F., Hosseini-Nasab, H., Mu{\~{n}}uzuri, J., \& Fakhrzad, M.~B. (2020).
\newblock {A multi-objective pharmaceutical supply chain network based on a robust fuzzy model: A comparison of meta-heuristics}.
\newblock {\em Applied Soft Computing Journal}, 92, 106331.
\newblock \url{https://doi.org/10.1016/j.asoc.2020.106331}.

\bibitem[Govindan et~al., 2017]{Govindan2017}
Govindan, K., Fattahi, M., \& Keyvanshokooh, E. (2017).
\newblock {Supply chain network design under uncertainty: A comprehensive review and future research directions}.
\newblock {\em European Journal of Operational Research}, 263(1), 108--141.
\newblock http://doi.org/10.1016/j.ejor.2017.04.009.

\bibitem[Hansen \& Grunow, 2015]{HANSEN2015}
Hansen, K. R.~N. \& Grunow, M. (2015).
\newblock Planning operations before market launch for balancing time-to-market and risks in pharmaceutical supply chains.
\newblock {\em International Journal of Production Economics}, 161, 129--139.
\newblock \url{https://doi.org/10.1016/j.ijpe.2014.10.010}.

\bibitem[Hasani \& Khosrojerdi, 2016]{Hasani2016}
Hasani, A. \& Khosrojerdi, A. (2016).
\newblock {Robust global supply chain network design under disruption and uncertainty considering resilience strategies: A parallel memetic algorithm for a real-life case study}.
\newblock {\em Transportation Research Part E: Logistics and Transportation Review}, 87, 20--52.
\newblock \url{https://doi.org/10.1016/j.tre.2015.12.009}.

\bibitem[Heese \& Kemahlioglu-Ziya, 2023]{Heese2023}
Heese, H.~S. \& Kemahlioglu-Ziya, E. (2023).
\newblock Capacity shortages, regulation, and firm incentives in the generic drugs industry.
\newblock {\em Naval Research Logistics}, 70(8), 811--831.
\newblock \url{https://doi.org/10.1002/nav.22132}.

\bibitem[Hilletofth et~al., 2019]{Hilletofth2019}
Hilletofth, P., Eriksson, D., Tate, W., \& Kinkel, S. (2019).
\newblock {Right-shoring: Making resilient offshoring and reshoring decisions}.
\newblock {\em Journal of Purchasing and Supply Management}, 25(3).
\newblock \url{https://doi.org/10.1016/j.pursup.2019.100540}.

\bibitem[Iacocca et~al., 2022]{Iacocca2022}
Iacocca, K., Mahar, S., \& {Daniel Wright}, P. (2022).
\newblock {Strategic horizontal integration for drug cost reduction in the pharmaceutical supply chain}.
\newblock {\em Omega}, 108, 17.
\newblock \url{https://doi.org/10.1016/j.omega.2021.102589}.

\bibitem[{IARC}, 2020]{IARC-WHO2020}
{IARC} (2020).
\newblock {Cancer over time}.
\newblock Retrieved from \url{https://gco.iarc.fr/overtime/en}. Accessed October 1, 2022.

\bibitem[Jahre et~al., 2016]{Jahre2016}
Jahre, M., Kembro, J., Rezvanian, T., Ergun, O., H{\aa}pnes, S.~J., \& Berling, P. (2016).
\newblock {Integrating supply chains for emergencies and ongoing operations in UNHCR}.
\newblock {\em Journal of Operations Management}, 45, 57--72.
\newblock \url{http://doi.org/10.1016/j.jom.2016.05.009}.

\bibitem[Kalish \& Wolf, 2021]{Kalish2021}
Kalish, I. \& Wolf, M. (2021).
\newblock {Supply chain resilience in the face of geopolitical risks}.
\newblock Retrieved from \url{https://www2.deloitte.com/us/en/insights/economy/us-china-trade-war-supply-chain.html}. Accessed May 5, 2023.

\bibitem[Kazaz et~al., 2023]{Kazaz2023}
Kazaz, B., Webster, S., \& Yadav, P. (2023).
\newblock Increasing the supply of health products in underserved regions.
\newblock {\em Production and Operations Management}, 32(12), 4212--4228.
\newblock \url{https://doi.org/10.1111/poms.14085}.

\bibitem[{Kchaou Boujelben} \& Boulaksil, 2018]{KchaouBoujelben2018}
{Kchaou Boujelben}, M. \& Boulaksil, Y. (2018).
\newblock {Modeling international facility location under uncertainty: A review, analysis, and insights}.
\newblock {\em IISE Transactions}, 50(6), 535--551.
\newblock \url{https://doi.org/10.1080/24725854.2017.1408165}.

\bibitem[Kleywegt et~al., 2002]{Kleywegt2002}
Kleywegt, A.~J., Shapiro, A., \& Homem-de Mello, T. (2002).
\newblock {The Sample Average Approximation Method for stochastic discrete optimization}.
\newblock {\em SIAM Journal on Optimization}, 12(2), 479--502.
\newblock \url{https://doi.org/10.1137/S1052623499363220}.

\bibitem[Kochan \& Nowicki, 2018]{Kochan2018}
Kochan, C.~G. \& Nowicki, D.~R. (2018).
\newblock {Supply chain resilience: A systematic literature review and typological framework}.
\newblock {\em International Journal of Physical Distribution and Logistics Management}, 48(8), 842--865.
\newblock \url{https://doi.org/10.1108/IJPDLM-02-2017-0099}.

\bibitem[Li \& Dong, 2022]{DONG-LI2022}
Li, D. \& Dong, C. (2022).
\newblock Government regulations to mitigate the shortage of life-saving goods in the face of a pandemic.
\newblock {\em European Journal of Operational Research}, 301(3), 942--955.
\newblock \url{https://doi.org/10.1016/j.ejor.2021.11.042}.

\bibitem[Li et~al., 2023a]{LiJinfeng2023}
Li, J., Liu, Y., \& Yang, G. (2023a).
\newblock {Two-stage distributionally robust optimization model for a pharmaceutical cold supply chain network design problem}.
\newblock {\em International Transactions in Operational Research}, 0(0), 1--35.
\newblock \url{https://doi.org/10.1111/itor.13267}.

\bibitem[Li et~al., 2020]{LiXiaopeng2020}
Li, X., Wang, H., Hao, G., \& Xia, C. (2020).
\newblock The mechanism of alliance promotes cooperation in the spatial multi-games.
\newblock {\em Physics Letters A}, 384(20), 126414.
\newblock \url{https://doi.org/10.1016/j.physleta.2020.126414}.

\bibitem[Li et~al., 2023b]{LiZhao2023}
Li, Z., Xia, T., Shen, W., \& Chen, S. (2023b).
\newblock {Research on co-opetition mechanism between pharmaceutical enterprises and third-party logistics in drug distribution of medical community}.
\newblock {\em International Journal of Environmental Research and Public Health}, 20(1), 19.
\newblock \url{https://doi.org/10.3390/ijerph20010609}.

\bibitem[Maharjan \& Kato, 2022]{Maharjan2022}
Maharjan, R. \& Kato, H. (2022).
\newblock {Resilient supply chain network design: A systematic literature review}.
\newblock {\em Transport Reviews}, 42(6), 739--761.
\newblock https://doi.org/10.1080/01441647.2022.2080773.

\bibitem[Marques et~al., 2020]{Marques2020}
Marques, C.~M., Moniz, S., de~Sousa, J.~P., Barbosa-Povoa, A.~P., \& Reklaitis, G. (2020).
\newblock {Decision-support challenges in the chemical-pharmaceutical industry: Findings and future research directions}.
\newblock {\em Computers and Chemical Engineering}, 134, 106672.
\newblock \url{https://doi.org/10.1016/j.compchemeng.2019.106672}.

\bibitem[Martei et~al., 2020]{Martei2020}
Martei, Y.~M., Iwamoto, K., Barr, R.~D., Wiernkowski, J.~T., \& Robertson, J. (2020).
\newblock {Shortages and price variability of essential cytotoxic medicines for treating children with cancers}.
\newblock {\em BMJ Global Health}, 5(11), e003282.
\newblock \url{https://doi.org/10.1136/bmjgh-2020-003282}.

\bibitem[McIvor \& Bals, 2021]{McIvor2021}
McIvor, R. \& Bals, L. (2021).
\newblock {A multi-theory framework for understanding the reshoring decision}.
\newblock {\em International Business Review}, 30(6), 101827.
\newblock https://doi.org/10.1016/j.ibusrev.2021.101827.

\bibitem[{McKinsey Global Institute}, 2020]{McKinseyGlobalInstitute2020}
{McKinsey Global Institute} (2020).
\newblock {Risk, resilience, and rebalancing in global value chains}.
\newblock Retrieved from \url{https://www.mckinsey.com/capabilities/operations/our-insights/risk-resilience-and-rebalancing-in-global-value-chains}. Accessed May 2, 2023.

\bibitem[Mohiuddin et~al., 2019]{Mohiuddin2019}
Mohiuddin, M., Rashid, M., Azad, S.~A., \& Su, Z. (2019).
\newblock {Back-shoring or re-shoring: Determinants of manufacturing offshoring from emerging to least developing countries}.
\newblock {\em International Journal Of Logistics: Research And Applications}, 22(1), 78--97.
\newblock \url{https://doi.org/10.1080/13675567.2018.1475554}.

\bibitem[Mulcahy et~al., 2021]{RR-2956-ASPEC}
Mulcahy, A.~W., Whaley, C.~M., Gizaw, M., Schwam, D., Edenfield, N., \& Becerra-Ornelas, A.~U. (2021).
\newblock {\em International Prescription Drug Price Comparisons: Current Empirical Estimates and Comparisons with Previous Studies}.
\newblock Santa Monica, CA: RAND Corporation.

\bibitem[Namdar et~al., 2018]{Namdar2018}
Namdar, J., Li, X., Sawhney, R., \& Pradhan, N. (2018).
\newblock {Supply chain resilience for single and multiple sourcing in the presence of disruption risks}.
\newblock {\em International Journal of Production Research}, 56(6), 2339--2360.
\newblock https://doi.org/10.1080/00207543.2017.1370149.

\bibitem[{NASEM}, 2022]{NASEM2022}
{NASEM} (2022).
\newblock {\em Building Resilience into the Nation's Medical Product Supply Chains}.
\newblock Washington, DC: The National Academies Press.
\newblock \url{https://doi.org/10.17226/26420}.

\bibitem[Nguyen et~al., 2021]{Nguyen2021}
Nguyen, H., Sharkey, T.~C., Wheeler, S., Mitchell, J.~E., \& Wallace, W.~A. (2021).
\newblock {Towards the development of quantitative resilience indices for multi-echelon assembly supply chains}.
\newblock {\em Omega}, 99, 102199.
\newblock https://doi.org/10.1016/j.omega.2020.102199.

\bibitem[Ouazine et~al., 2018]{Ouazine2018}
Ouazine, K., Slimani, H., \& Tari, A. (2018).
\newblock {Alliances in graphs: Parameters, properties and applications—A survey}.
\newblock {\em AKCE International Journal of Graphs and Combinatorics}, 15(2), 115--154.
\newblock \url{https://doi.org/10.1016/j.akcej.2017.05.002}.

\bibitem[Rawls \& Turnquist, 2010]{Rawls2010}
Rawls, C.~G. \& Turnquist, M.~A. (2010).
\newblock {Pre-positioning of emergency supplies for disaster response}.
\newblock {\em Transportation Research Part B: Methodological}, 44(4), 521--534.
\newblock http://doi.org/10.1016/j.trb.2009.08.003.

\bibitem[{Sabogal De La Pava} \& Tucker, 2022]{SabogalDeLaPava2022}
{Sabogal De La Pava}, M.~L. \& Tucker, E.~L. (2022).
\newblock {Drug shortages in low- and middle-income countries: Colombia as a case study}.
\newblock {\em Journal of Pharmaceutical Policy and Practice}, 15(1), 1--8.
\newblock \url{https://doi.org/10.1186/s40545-022-00439-7}.

\bibitem[Sabouhi et~al., 2018]{Sabouhi2018}
Sabouhi, F., Pishvaee, M.~S., \& Jabalameli, M.~S. (2018).
\newblock {Resilient supply chain design under operational and disruption risks considering quantity discount: A case study of pharmaceutical supply chain}.
\newblock {\em Computers and Industrial Engineering}, 126, 657--672.
\newblock \url{https://doi.org/10.1016/j.cie.2018.10.001}.

\bibitem[Saedi et~al., 2016]{SAEDI2016}
Saedi, S., Kundakcioglu, O.~E., \& Henry, A.~C. (2016).
\newblock Mitigating the impact of drug shortages for a healthcare facility: An inventory management approach.
\newblock {\em European Journal of Operational Research}, 251(1), 107--123.
\newblock \url{https://doi.org/10.1016/j.ejor.2015.11.017}.

\bibitem[Santoso et~al., 2005]{SANTOSO200596}
Santoso, T., Ahmed, S., Goetschalckx, M., \& Shapiro, A. (2005).
\newblock A stochastic programming approach for supply chain network design under uncertainty.
\newblock {\em European Journal of Operational Research}, 167(1), 96--115.
\newblock \url{https://doi.org/10.1016/j.ejor.2004.01.046}.

\bibitem[Simon, 2022]{Simon2022}
Simon, D.~W. (2022).
\newblock {Managing supply chain disruption in an era of geopolitical risk}.
\newblock Retrieved from \url{https://www.foley.com/en/insights/publications/2022/07/avoid-supply-chain-disruption-geopolitical-risk}. Accessed May 2, 2023.

\bibitem[Singh et~al., 2019]{Singh2019}
Singh, C.~S., Soni, G., \& Badhotiya, G.~K. (2019).
\newblock {Performance indicators for supply chain resilience: Review and conceptual framework}.
\newblock {\em Journal of Industrial Engineering International}, 15, 105--117.
\newblock http://doi.org/10.1007/s40092-019-00322-2.

\bibitem[Singh et~al., 2018]{Singh2018}
Singh, H., Garg, R.~K., \& Sachdeva, A. (2018).
\newblock {Supply chain collaboration: A state-of-the-art literature review}.
\newblock {\em Uncertain Supply Chain Management}, 6, 149--180.
\newblock \url{https://doi.org/10.5267/j.uscm.2017.8.002}.

\bibitem[Snyder et~al., 2016]{Snyder2016}
Snyder, L.~V., Atan, Z., Peng, P., Rong, Y., Schmitt, A.~J., \& Sinsoysal, B. (2016).
\newblock {OR/MS models for supply chain disruptions: A review}.
\newblock {\em IIE Transactions}, 48(2), 89--109.
\newblock https://doi.org/10.1080/0740817X.2015.1067735.

\bibitem[Sousa et~al., 2011]{Sousa2011}
Sousa, R.~T., Liu, S., Papageorgiou, L.~G., \& Shah, N. (2011).
\newblock {Global supply chain planning for pharmaceuticals}.
\newblock {\em Chemical Engineering Research and Design}, 89(11), 2396--2409.
\newblock \url{https://doi.org/10.1016/j.cherd.2011.04.005}.

\bibitem[Suda et~al., 2022]{Suda2022}
Suda, K.~J., Kim, K.~C., Hernandez, I., Gellad, W.~F., Rothenberger, S., Campbell, A., Malliart, L., \& Tadrous, M. (2022).
\newblock {The global impact of COVID-19 on drug purchases: A cross-sectional time series analysis}.
\newblock {\em Journal of the American Pharmacists Association}, 62(3), 766--774.
\newblock \url{https://doi.org/10.1016/j.japh.2021.12.014}.

\bibitem[Suryawanshi \& Dutta, 2022]{Suryawanshi2022}
Suryawanshi, P. \& Dutta, P. (2022).
\newblock {Optimization models for supply chains under risk, uncertainty, and resilience: A state-of-the-art review and future research directions}.
\newblock {\em Transportation Research Part E: Logistics and Transportation Review}, 157(Dec 2021), 1--45.
\newblock \url{https://doi.org/10.1016/j.tre.2021.102553}.

\bibitem[Susarla \& Karimi, 2012]{Susarla2012}
Susarla, N. \& Karimi, I.~A. (2012).
\newblock {Integrated supply chain planning for multinational pharmaceutical enterprises}.
\newblock {\em Computers and Chemical Engineering}, 42, 168--177.
\newblock \url{https://doi.org/10.1016/j.compchemeng.2012.03.002}.

\bibitem[Tatar et~al., 2022]{Tatar2022}
Tatar, M., Shoorekchali, J.~M., Faraji, M.~R., Seyyedkolaee, M.~A., Pag{\'{a}}n, J.~A., \& Wilson, F.~A. (2022).
\newblock {COVID-19 vaccine inequality: A global perspective}.
\newblock {\em Journal of Global Health}, 12, 10--13.
\newblock \url{https://doi.org/10.7189/jogh.12.03072}.

\bibitem[{The White House}, 2021]{TheWhiteHouse2021}
{The White House} (2021).
\newblock {\em {Building Resilient Supply Chains, Revitalizing American Manufacturing, and Fostering Broad-Based Growth: 100-Day Reviews Under Executive Order 14017}}.
\newblock Washington, DC: {The White House}.

\bibitem[Tordecilla et~al., 2021]{Tordecilla2021}
Tordecilla, R.~D., Juan, A.~A., Montoya-Torres, J.~R., Quintero-Araujo, C.~L., \& Panadero, J. (2021).
\newblock {Simulation-optimization methods for designing and assessing resilient supply chain networks under uncertainty scenarios: A review}.
\newblock {\em Simulation Modelling Practice and Theory}, 106(Aug 2020), 102166.
\newblock http://doi.org/10.1016/j.simpat.2020.102166.

\bibitem[{Trading Economics}, 2022]{TradingEconomics2022}
{Trading Economics} (2022).
\newblock {Capacity utilization}.
\newblock Retrieved from \url{https://tradingeconomics.com/country-list/capacity-utilization}. Accessed December 10, 2022.

\bibitem[Tucker \& Daskin, 2022]{Tucker2022}
Tucker, E.~L. \& Daskin, M.~S. (2022).
\newblock Pharmaceutical supply chain reliability and effects on drug shortages.
\newblock {\em Computers and Industrial Engineering}, 169, 108258.
\newblock \url{https://doi.org/10.1016/j.cie.2022.108258}.

\bibitem[Tucker et~al., 2020]{Tucker2020}
Tucker, E.~L., Daskin, M.~S., Sweet, B.~V., \& Hopp, W.~J. (2020).
\newblock {Incentivizing resilient supply chain design to prevent drug shortages: Policy analysis using two- and multi-stage stochastic programs}.
\newblock {\em IISE Transactions}, 52(4), 394--412.
\newblock \url{https://doi.org/10.1080/24725854.2019.1646441}.

\bibitem[{UN Comtrade Database}, 2019]{UNComtradeDatabase2019}
{UN Comtrade Database} (2019).
\newblock {Global trade data}.
\newblock Retrieved from \url{https://comtradeplus.un.org/}. Accessed December 10, 2022.

\bibitem[{UNCTADstat}, 2016]{UNCTADstat2016}
{UNCTADstat} (2016).
\newblock {International merchandise trade - transport costs}.
\newblock Retrieved from \url{https://unctadstat.unctad.org/wds/ReportFolders/reportFolders.aspx?sCS\_ChosenLang=en}. Accessed February 17, 2023.

\bibitem[{UUDIS}, 2022]{Utah2022}
{UUDIS} (2022).
\newblock {Current Drug Shortages}.
\newblock Retrieved from \url{https://www.ashp.org/drug-shortages/current-shortages}. Accessed May 3, 2023.

\bibitem[Varadarajan \& Cunningham, 1995]{Varadarajan1995}
Varadarajan, P.~R. \& Cunningham, M.~H. (1995).
\newblock {Strategic alliances: A synthesis of conceptual foundations}.
\newblock {\em Journal of the Academy of Marketing Science}, 23, 282--296.
\newblock \url{https://doi.org/10.1177/009207039502300408}.

\bibitem[Verweij et~al., 2003]{Verweij2003}
Verweij, B., Ahmed, S., Kleywegt, A.~J., Nemhauser, G., \& Shapiro, A. (2003).
\newblock {The Sample Average Approximation Method applied to stochastic routing problems: A computational study}.
\newblock {\em Computational Optimization and Applications}, 24(2-3), 289--333.
\newblock \url{https://doi.org/10.1023/A:1021814225969}.

\bibitem[{World Bank}, 2020]{WorldBank2020}
{World Bank} (2020).
\newblock {\em {Doing Business 2020}}.
\newblock Washington, DC: {World Bank}.

\bibitem[WTO, 2020]{WTO2020-ER}
WTO (2020).
\newblock {Export prohibitions and restrictions}.
\newblock Retrieved from \url{https://doi.org/10.30875/827540bd-en}. Accessed May 5, 2023.

\bibitem[{WTO}, 2022a]{WTO2022a}
{WTO} (2022a).
\newblock {Quantitative restrictions database}.
\newblock Retrieved from \url{https://qr.wto.org/en\#/home}. Accessed October 1, 2022.

\bibitem[{WTO}, 2022b]{WTO2022b}
{WTO} (2022b).
\newblock {Trade monitoring database}.
\newblock Retrieved from \url{https://tmdb.wto.org/en/explore/goods\#nogo\%0A}. Accessed October 1, 2022.

\bibitem[WTO, 2024]{WTO2024}
WTO (2024).
\newblock {Regional Trade Agreements}.
\newblock Retrieved from \url{https://rtais.wto.org/UI/PublicMaintainRTAHome.aspx}. Accessed August 3, 2024.

\bibitem[Zandkarimkhani et~al., 2020]{Zandkarimkhani2020}
Zandkarimkhani, S., Mina, H., Biuki, M., \& Govindan, K. (2020).
\newblock {A chance constrained fuzzy goal programming approach for perishable pharmaceutical supply chain network design}.
\newblock {\em Annals of Operations Research}, 295(1), 425--452.
\newblock \url{https://doi.org/10.1007/s10479-020-03677-7}.

\bibitem[Zhao, 2023]{Zhao-tutorial}
Zhao, H. (2023).
\newblock Pharmaceutical supply chains and drug shortages.
\newblock {\em INFORMS TutORials in Operations Research}, 0(0), 228--245.
\newblock \url{https://doi.org/10.1287/educ.2023.0258}.

\end{thebibliography}
